\documentclass{commat}

\DeclareMathOperator{\id}{id}
\newcommand{\idbar}{\mathop{\overline{\id}}}

\title{On the square of the antipode in a connected filtered Hopf algebra}

\author{Darij Grinberg}

\affiliation{
    \address{%
    Darij Grinberg --
    Drexel University, Department of Mathematics,
    Korman Center 263, 15 S 33rd Street,
    Philadelphia, PA 19104, USA
    }
    \email{%
    darijgrinberg@gmail.com
    }
    }

\abstract{%
It is well-known that the antipode $S$ of a commutative or
cocommutative Hopf algebra satisfies $S^{2}=\operatorname*{id}$ (where
$S^{2}=S\circ S$). Recently, similar results have been obtained by Aguiar,
Lauve and Mahajan for connected graded Hopf algebras: Namely, if $H$ is a
connected graded Hopf algebra with grading $H=\bigoplus_{n\geq0}H_{n}$, then
each positive integer $n$ satisfies $\left(  \operatorname*{id}-S^{2}\right)
^{n}\left(  H_{n}\right)  =0$ and (even stronger)
\[
\left(  \left(  \operatorname*{id}+S\right)  \circ\left(  \operatorname*{id}-S^{2}\right)^{n-1}\right)  \left(  H_{n}\right) = 0.
\]
For some specific $H$'s such as the
Malvenuto--Reutenauer Hopf algebra $\operatorname*{FQSym}$, the exponents can
be lowered.

In this note, we generalize these results in several directions: We replace
the base field by a commutative ring, replace the Hopf algebra by a coalgebra
(actually, a slightly more general object, with no coassociativity required),
and replace both $\operatorname*{id}$ and $S^{2}$ by \textquotedblleft
coalgebra homomorphisms\textquotedblright\ (of sorts). Specializing back to
connected graded Hopf algebras, we show that the exponent $n$ in the identity
$\left(  \operatorname*{id}-S^{2}\right)  ^{n}\left(  H_{n}\right)  =0$ can be
lowered to $n-1$ (for $n>1$) if and only if $\left(  \operatorname*{id}%
-S^{2}\right)  \left(  H_{2}\right)  =0$. (A sufficient condition for this is
that every pair of elements of $H_{1}$ commutes; this is satisfied, e.g., for
$\operatorname*{FQSym}$.)
}



\VOLUME{31}
\YEAR{2023}
\NUMBER{1}
\firstpage{303}
\DOI{https://doi.org/10.46298/cm.10431}

\begin{paper}

Consider, for simplicity, a connected graded Hopf algebra $H$ over a field (we
will soon switch to more general settings). Let $S$ be the antipode of $H$. A
classical result (e.g., \cite[Proposition 4.0.1 6)]{Sweedler} or
\cite[Corollary 3.3.11]{HaGuKi10} or \cite[Theorem 2.1.4 (vi)]{Abe80} or
\cite[Corollary 7.1.11]{Radfor12}) says that $S^{2}=\operatorname*{id}$ when
$H$ is commutative or cocommutative. (Here and in the following, powers are
composition powers; thus, $S^{2}$ means $S\circ S$.) In general,
$S^{2}=\operatorname*{id}$ need not hold. However, in \cite[Proposition
7]{AguLau14}, Aguiar and Lauve showed that $S^{2}$ is still locally unipotent,
and more precisely, we have
\[
\left(  \operatorname*{id}-S^{2}\right)  ^{n}\left(  H_{n}\right)
=0\ \ \ \ \ \ \ \ \ \ \text{for each }n>0,
\]
where $H_{n}$ denotes the $n$-th graded component of $H$. Later, Aguiar and
Mahajan \cite[Lemma~12.50]{Aguiar17} strengthened this equality to
\[
\left(  \left(  \operatorname*{id}+S\right)  \circ\left(  \operatorname*{id}%
-S^{2}\right)  ^{n-1}\right)  \left(  H_{n}\right)
=0\ \ \ \ \ \ \ \ \ \ \text{for each }n>0.
\]
For specific combinatorially interesting Hopf algebras, even stronger results
hold; in particular,
\[
\left(  \operatorname*{id}-S^{2}\right)  ^{n-1}\left(  H_{n}\right)
=0\ \ \ \ \ \ \ \ \ \ \text{holds for each }n>1
\]
when $H$ is the Malvenuto--Reutenauer Hopf algebra (see \cite[Example~8]{AguLau14}).

In this note, we will unify these results and transport them to a much more
general setting. First of all, the ground field will be replaced by an
arbitrary commutative ring; this generalization is not unexpected, but renders
the proof strategy of \cite[Proposition~7]{AguLau14} insufficient\footnote{In
fact, the proof in \cite[Proposition~7]{AguLau14} relies on the coradical
filtration of $H$ and its associated graded structure $\operatorname*{gr}H$.
If the base ring is a field, then $\operatorname*{gr}H$ is a well-defined
commutative Hopf algebra (see, e.g., \cite[Lemma~1]{AguLau14}), and thus the
antipode of $H$ can be viewed as a \textquotedblleft
deformation\textquotedblright\ of the antipode of $\operatorname*{gr}H$. But
the latter antipode does square to $\operatorname*{id}$ because
$\operatorname*{gr}H$ is commutative. Unfortunately, this proof does not
survive our generalization; in fact, even defining a Hopf algebra structure on
$\operatorname*{gr}H$ would likely require at least some flatness
assumptions.}. Second, we will replace the Hopf algebra by a coalgebra, or
rather by a more general structure that does not even require coassociativity.
The squared antipode $S^{2}$ will be replaced by an arbitrary
\textquotedblleft coalgebra\textquotedblright\ endomorphism\ $f$ (we are using
scare quotes because our structure is not really a coalgebra), and the
identity map by another such endomorphism\ $e$. Finally, the graded components
will be replaced by an arbitrary sequence of submodules satisfying certain
compatibility relations. We state the general result in Section
\ref{sec.general} and prove it in Section \ref{sec.pf.thm.id-f.gen}. In
Sections \ref{sec.cfc}--\ref{sec.hopf}, we progressively specialize this
result: first to connected filtered coalgebras with coalgebra endomorphisms
(in Section \ref{sec.cfc}), then to connected filtered Hopf algebras with
$S^{2}$ (in Section \ref{sec.filhopf}), and finally to connected graded Hopf
algebras with $S^{2}$ (in Section \ref{sec.hopf}). The latter specialization
covers the results of Aguiar and Lauve. (The Malvenuto--Reutenauer Hopf
algebra turns out to be a red herring; any connected graded Hopf algebra $H$
with the property that $ab=ba$ for all $a,b\in H_{1}$ will do.)



\section{Notations}

We will use the notions of coalgebras, bialgebras and Hopf algebras over a
commutative ring, as defined (e.g.) in \cite[Chapter 2]{Abe80}, \cite[Chapter
1]{GriRei}, \cite[Chapters 2, 3]{HaGuKi10}, \cite[Chapters 2, 5, 7]{Radfor12}
or \cite[Chapters I--IV]{Sweedler}. (In particular, our Hopf algebras are
\textbf{not} twisted by a $\mathbb{Z}/2$-grading as the topologists' ones
are.) We use the same notations for Hopf algebras as in \cite[Chapter
1]{GriRei}. In particular:

\begin{itemize}
\item We let $\mathbb{N}=\left\{  0,1,2,\ldots\right\}  $.

\item \textquotedblleft Rings\textquotedblright\ and \textquotedblleft
algebras\textquotedblright\ are always required to be associative and have a unity.

\item We fix a commutative ring $\mathbf{k}$. The symbols \textquotedblleft%
$\otimes$\textquotedblright\ and \textquotedblleft$\operatorname*{End}%
$\textquotedblright\ shall always mean \textquotedblleft$\otimes_{\mathbf{k}}%
$\textquotedblright\ and \textquotedblleft$\operatorname*{End}%
\nolimits_{\mathbf{k}}$\textquotedblright, respectively. The unity of the ring
$\mathbf{k}$ will be called $1_{\mathbf{k}}$ or just $1$ if confusion is unlikely.

\item The comultiplication and the counit of a $\mathbf{k}$-coalgebra are
denoted by $\Delta$ and $\epsilon$.

\item \textquotedblleft Graded\textquotedblright\ $\mathbf{k}$-modules mean
$\mathbb{N}$-graded $\mathbf{k}$-modules. The base ring $\mathbf{k}$ itself is
not supposed to have any nontrivial grading.

\item The $n$-th graded component of a graded $\mathbf{k}$-module $V$ will be
called $V_{n}$. If $n<0$, then this is the zero submodule $0$.

\item A graded $\mathbf{k}$-Hopf algebra means a $\mathbf{k}$-Hopf algebra
that has a grading as a $\mathbf{k}$-module, and whose structure maps
(multiplication, unit, comultiplication and counit) are graded maps. (The
antipode is automatically graded in this case, by \cite[Exercise 1.4.29
(e)]{GriRei}.)

\item If $f$ is a map from a set to itself, and if $k\in\mathbb{N}$ is
arbitrary, then $f^{k}$ shall denote the map $\underbrace{f\circ f\circ
\cdots\circ f}_{k\text{ times}}$. (Thus, $f^{1}=f$ and $f^{0}%
=\operatorname*{id}$.)
\end{itemize}

\section{Theorems}

\subsection{\label{sec.general}The main theorem}

We can now state the main result of this note:

\begin{theorem}
\label{thm.id-f.gen}Let $D$ be a $\mathbf{k}$-module, and let $\left(
D_{1},D_{2},D_{3},\ldots\right)  $ be a sequence of $\mathbf{k}$-submodules of
$D$. Let $\delta:D\rightarrow D\otimes D$ be any $\mathbf{k}$-linear map.

Let $e:D\rightarrow D$ and $f:D\rightarrow D$ be two $\mathbf{k}$-linear maps
such that%
\begin{align}
\operatorname*{Ker}\delta &  \subseteq\operatorname*{Ker}\left(  e-f\right)
\ \ \ \ \ \ \ \ \ \ \text{and}\label{pf.thm.id-f.gen.ass-Ker}\\
\left(  f\otimes f\right)  \circ\delta &  =\delta\circ
f\ \ \ \ \ \ \ \ \ \ \ \text{and}\label{pf.thm.id-f.gen.ass-mor}\\
\left(  e\otimes e\right)  \circ\delta &  =\delta\circ
e\ \ \ \ \ \ \ \ \ \ \text{and}\label{pf.thm.id-f.gen.ass-mor'}\\
f\circ e  &  =e\circ f. \label{pf.thm.id-f.gen.ass-comm}%
\end{align}

Let $p$ be a positive integer such that%
\begin{equation}
\left(  e-f\right)  \left(  D_{1}+D_{2}+\cdots+D_{p}\right)  =0.
\label{pf.thm.id-f.gen.ass-ann}%
\end{equation}
Assume furthermore that%
\begin{equation}
\delta\left(  D_{n}\right)  \subseteq\sum_{i=1}^{n-1}D_{i}\otimes
D_{n-i}\ \ \ \ \ \ \ \ \ \ \text{for each }n>p. \label{pf.thm.id-f.gen.ass-gr}%
\end{equation}
(Here, the \textquotedblleft$D_{i}\otimes D_{n-i}$\textquotedblright\ on the
right hand side means the image of $D_{i}\otimes D_{n-i}$ under the canonical
map $D_{i}\otimes D_{n-i}\rightarrow D\otimes D$ that is obtained by tensoring
the two inclusion maps $D_{i}\rightarrow D$ and $D_{n-i}\rightarrow D$
together. When $\mathbf{k}$ is not a field, this canonical map may fail to be injective.)

Then, for any integer $u>p$, we have%
\begin{equation}
\left(  e-f\right)  ^{u-p}\left(  D_{u}\right)  \subseteq\operatorname*{Ker}%
\delta\label{eq.thm.id-f.gen.clm1}%
\end{equation}
and%
\begin{equation}
\left(  e-f\right)  ^{u-p+1}\left(  D_{u}\right)  =0.
\label{eq.thm.id-f.gen.clm2}%
\end{equation}

\end{theorem}

As the statement of this theorem is not very intuitive, some explanations are
in order. The reader may think of the $D$ in Theorem \ref{thm.id-f.gen} as a
\textquotedblleft pre-coalgebra\textquotedblright, with $\delta$ being its
\textquotedblleft reduced coproduct\textquotedblright. Indeed, the easiest way
to obtain a nontrivial example is to fix a connected graded Hopf algebra $H$,
then define $D$ to be either $H$ or the \textquotedblleft positive
part\textquotedblright\ of $H$ (that is, the submodule $\bigoplus_{n>0}H_{n}$
of $H$), and define $\delta$ to be the map 
\[
    x\mapsto\Delta\left(  x\right) -x\otimes1-1\otimes x+\epsilon\left(  x\right)  1\otimes1
\]
(the so-called \emph{reduced coproduct} of $H$). From this point of view,
$\operatorname*{Ker}\delta$ can be regarded as the set of \textquotedblleft
primitive\textquotedblright\ elements of $D$. The maps $f$ and $e$ can be
viewed as two commuting \textquotedblleft coalgebra
endomorphisms\textquotedblright\ of $D$ (indeed, the assumptions
(\ref{pf.thm.id-f.gen.ass-mor}) and (\ref{pf.thm.id-f.gen.ass-mor'}) are
essentially saying that $f$ and $e$ preserve the \textquotedblleft reduced
coproduct\textquotedblright\ $\delta$). The submodules $D_{1},D_{2}%
,D_{3},\ldots$ are an analogue of the (positive-degree) graded components of
$D$, while the assumption (\ref{pf.thm.id-f.gen.ass-gr}) says that the
\textquotedblleft reduced coproduct\textquotedblright\ $\delta$
\textquotedblleft respects the grading\textquotedblright\ (as is indeed the
case for connected graded Hopf algebras).

We stress that $p$ is allowed to be $1$ in Theorem \ref{thm.id-f.gen}; in this
case, the assumption (\ref{pf.thm.id-f.gen.ass-ann}) simplifies to $\left(
e-f\right)  \left(  0\right)  =0$, which is automatically true by the
linearity of $e-f$.

We shall prove Theorem \ref{thm.id-f.gen} in Section \ref{sec.pf.thm.id-f.gen}%
. First, however, let us explore its consequences for coalgebras and Hopf
algebras, recovering in particular the results of Aguiar and Lauve promised in
the introduction.

\subsection{\label{sec.cfc}Connected filtered coalgebras}

We begin by specializing Theorem \ref{thm.id-f.gen} to the setting of a
connected filtered coalgebra. There are several ways to define what a filtered
coalgebra is; ours is probably the most liberal:

\begin{definition}
\label{def.fil-coal}A \emph{filtered }$\mathbf{k}$\emph{-coalgebra} means a
$\mathbf{k}$-coalgebra $C$ equipped with an infinite sequence $\left(
C_{\leq0},C_{\leq1},C_{\leq2},\ldots\right)  $ of $\mathbf{k}$-submodules of
$C$ satisfying the following three conditions:

\begin{itemize}
\item We have%
\begin{equation}
C_{\leq0}\subseteq C_{\leq1}\subseteq C_{\leq2}\subseteq\cdots.
\label{eq.def.fil-coal.chain}%
\end{equation}

\item We have%
\begin{equation}
\bigcup_{n\in\mathbb{N}}C_{\leq n}=C. \label{eq.def.fil-coal.union}%
\end{equation}

\item We have%
\begin{equation}
\Delta\left(  C_{\leq n}\right)  \subseteq\sum_{i=0}^{n}C_{\leq i}\otimes
C_{\leq n-i}\ \ \ \ \ \ \ \ \ \ \text{for each }n\in\mathbb{N}.
\label{eq.def.fil-coal.Del}%
\end{equation}
(Here, the \textquotedblleft$C_{\leq i}\otimes C_{\leq n-i}$\textquotedblright%
\ on the right hand side means the image of $C_{\leq i}\otimes C_{\leq n-i}$
under the canonical map $C_{\leq i}\otimes C_{\leq n-i}\rightarrow C\otimes C$
that is obtained by tensoring the two inclusion maps $C_{\leq i}\rightarrow C$
and $C_{\leq n-i}\rightarrow C$ together. When $\mathbf{k}$ is not a field,
this canonical map may fail to be injective.)
\end{itemize}

The sequence $\left(  C_{\leq0},C_{\leq1},C_{\leq2},\ldots\right)  $ is called
the \emph{filtration} of the filtered $\mathbf{k}$-coalgebra $C$.
\end{definition}

A more categorically-minded person might replace the condition~\eqref{eq.def.fil-coal.Del} by a stronger requirement (e.g., asking $\Delta$ to factor
through a linear map $C_{\leq n}\rightarrow\bigoplus_{i=0}^{n}C_{\leq
i}\otimes C_{\leq n-i}$, where the \textquotedblleft$\otimes$%
\textquotedblright\ signs now signify the actual tensor products rather than
their images in $C\otimes C$). However, we have no need for such stronger
requirements. Mercifully, all reasonable definitions of filtered $\mathbf{k}%
$-coalgebras agree when $\mathbf{k}$ is a field.

The condition (\ref{eq.def.fil-coal.union}) in Definition \ref{def.fil-coal}
shall never be used in the following; we merely state it to avoid muddling the
meaning of \textquotedblleft filtered $\mathbf{k}$-coalgebra\textquotedblright.

A graded $\mathbf{k}$-coalgebra $C$ automatically becomes a filtered
$\mathbf{k}$-coalgebra; indeed, we can define its filtration $\left(
C_{\leq0},C_{\leq1},C_{\leq2},\ldots\right)  $ by setting%
\[
C_{\leq n}=\bigoplus_{i=0}^{n}C_{i}\ \ \ \ \ \ \ \ \ \ \text{for all }%
n\in\mathbb{N}.
\]

\begin{definition}
\label{def.con-fil-coal}Let $C$ be a filtered $\mathbf{k}$-coalgebra with
filtration $\left(  C_{\leq0},C_{\leq1},C_{\leq2},\ldots\right)  $. Let
$1_{\mathbf{k}}$ denote the unity of the ring $\mathbf{k}$.

\textbf{(a)} The filtered $\mathbf{k}$-coalgebra $C$ is said to be
\emph{connected} if the restriction $\epsilon\mid_{C_{\leq0}}$ is a
$\mathbf{k}$-module isomorphism from $C_{\leq0}$ to $\mathbf{k}$.

\textbf{(b)} In this case, the element $\left(  \epsilon\mid_{C_{\leq0}%
}\right)  ^{-1}\left(  1_{\mathbf{k}}\right)  \in C_{\leq0}$ is called the
\emph{unity} of $C$ and is denoted by $1_{C}$.

Now, assume that $C$ is a connected filtered $\mathbf{k}$-coalgebra.

\textbf{(c)} An element $x$ of $C$ is said to be \emph{primitive} if
$\Delta\left(  x\right)  =x\otimes1_{C}+1_{C}\otimes x$.

\textbf{(d)} The set of all primitive elements of $C$ is denoted by
$\operatorname*{Prim}C$.
\end{definition}

These notions of \textquotedblleft connected\textquotedblright,
\textquotedblleft unity\textquotedblright\ and \textquotedblleft
primitive\textquotedblright\ specialize to the commonly established concepts
of these names when $C$ is a graded $\mathbf{k}$-bialgebra. Indeed, Definition
\ref{def.con-fil-coal} \textbf{(b)} defines the unity $1_{C}$ of $C$ to be the
unique element of $C_{\leq0}$ that gets sent to $1_{\mathbf{k}}$ by the map
$\epsilon$; but this property is satisfied for the unity of a graded
$\mathbf{k}$-bialgebra as well. (We will repeat this argument in more detail
later on, in the proof of Proposition \ref{prop.fil-bial.1=1}.)

The following property of connected filtered $\mathbf{k}$-coalgebras will be
crucial for us:

\begin{proposition}
\label{prop.cfc.delta2}Let $C$ be a connected filtered $\mathbf{k}$-coalgebra,
and let $\left(  C_{\leq0},C_{\leq1},C_{\leq2},\ldots\right)  $ be its filtration.
Define a $\mathbf{k}$-linear map $\delta:C\rightarrow C\otimes C$ by setting
\[
\delta\left(  c\right)  :=\Delta\left(  c\right)  -c\otimes1_{C}-1_{C}\otimes
c+\epsilon\left(  c\right)  1_{C}\otimes1_{C}\ \ \ \ \ \ \ \ \ \ \text{for
each }c\in C.
\]
Then:

\textbf{(a)} We have
\[
\delta\left(  C_{\leq n}\right)  \subseteq\sum_{i=1}^{n-1}C_{\leq i}\otimes
C_{\leq n-i}\ \ \ \ \ \ \ \ \ \ \text{for each }n>0.
\]

\textbf{(b)} If $f:C\rightarrow C$ is a $\mathbf{k}$-coalgebra homomorphism
satisfying $f\left(  1_{C}\right)  =1_{C}$, then we have $\left(  f\otimes
f\right)  \circ\delta=\delta\circ f$.

\textbf{(c)} We have $\operatorname*{Prim}C=\left(  \operatorname*{Ker}%
\delta\right)  \cap\left(  \operatorname*{Ker}\epsilon\right)  $.

\textbf{(d)} The set $\operatorname*{Prim}C$ is a $\mathbf{k}$-submodule of
$C$.

\textbf{(e)} We have $\operatorname*{Ker}\delta=\mathbf{k}\cdot1_{C}%
+\operatorname*{Prim}C$.
\end{proposition}

The map $\delta$ in Proposition \ref{prop.cfc.delta2} is called the
\emph{reduced coproduct} of $C$. Nothing in Proposition \ref{prop.cfc.delta2}
is really novel; in particular, results similar to Proposition
\ref{prop.cfc.delta2} \textbf{(a)} have appeared all over the literature
(e.g., in \cite[Lemma 1.3.6 (1)]{HecSch20} with the right hand side
$\sum_{i=1}^{n-1}C_{\leq i}\otimes C_{\leq n-i}$ replaced by the weaker
$C_{\leq n-1}\otimes C_{\leq n-1}$, or in \cite[Lemma 1.3.6 (2)]{HecSch20} for
graded coalgebras, or in \cite[Proposition 10.0.2]{Sweedler} for a specific
filtration defined over a field, or in \cite[(3.2.6)]{Mancho06} for
$\delta\left(  C_{\leq n}\cap\operatorname*{Ker}\varepsilon\right)  $ instead
of $\delta\left(  C_{\leq n}\right)  $), and the arguments used in these
sources can often be repurposed with some care to apply to Proposition
\ref{prop.cfc.delta2} \textbf{(a)}. Parts \textbf{(b)}--\textbf{(e)} of the
proposition are folklore as well. For the sake of completeness, we shall
nevertheless prove Proposition \ref{prop.cfc.delta2} in Section
\ref{sec.pf.prop.cfc.delta2}, which any reader with subject experience can skip.

Proposition \ref{prop.cfc.delta2} helps us apply Theorem \ref{thm.id-f.gen} to
filtered $\mathbf{k}$-coalgebras, resulting in the following:

\begin{corollary}
\label{cor.id-f.cfc}Let $C$ be a connected filtered $\mathbf{k}$-coalgebra
with filtration $\left(  C_{\leq0},C_{\leq1},C_{\leq2},\ldots\right)  $.

Let $e:C\rightarrow C$ and $f:C\rightarrow C$ be two $\mathbf{k}$-coalgebra
homomorphisms such that%
\begin{align}
e\left(  1_{C}\right)   &  =1_{C}\ \ \ \ \ \ \ \ \ \ \text{and}\nonumber\\
f\left(  1_{C}\right)   &  =1_{C}\ \ \ \ \ \ \ \ \ \ \text{and}\nonumber\\
\operatorname*{Prim}C  &  \subseteq\operatorname*{Ker}\left(  e-f\right)
\ \ \ \ \ \ \ \ \ \ \text{and}\label{eq.cor.id-f.cfc.ass-Ker}\\
f\circ e  &  =e\circ f. \label{eq.cor.id-f.cfc.ass-comm}%
\end{align}

Let $p$ be a positive integer such that%
\begin{equation}
\left(  e-f\right)  \left(  C_{\leq p}\right)  =0.
\label{eq.cor.id-f.cfc.ass-ann}%
\end{equation}

Then:

\textbf{(a)} For any integer $u>p$, we have%
\begin{equation}
\left(  e-f\right)  ^{u-p}\left(  C_{\leq u}\right)  \subseteq
\operatorname*{Prim}C. \label{eq.cor.id-f.cfc.claim-1}%
\end{equation}

\textbf{(b)} For any integer $u\geq p$, we have%
\begin{equation}
\left(  e-f\right)  ^{u-p+1}\left(  C_{\leq u}\right)  =0.
\label{eq.cor.id-f.cfc.claim-2}%
\end{equation}

\end{corollary}

Corollary \ref{cor.id-f.cfc} results from an easy (although not completely
immediate) application of Theorem \ref{thm.id-f.gen} and Proposition
\ref{prop.cfc.delta2}. A detailed proof can be found in Section
\ref{sec.pf.cor.id-f.cfc}.

Specializing Corollary \ref{cor.id-f.cfc} further to the case of $p=1$, we can
obtain a nicer result:

\begin{corollary}
\label{cor.id-f.cfc1}Let $C$ be a connected filtered $\mathbf{k}$-coalgebra
with filtration $\left(  C_{\leq0},C_{\leq1},C_{\leq2},\ldots\right)  $.

Let $e:C\rightarrow C$ and $f:C\rightarrow C$ be two $\mathbf{k}$-coalgebra
homomorphisms such that%
\begin{align*}
e\left(  1_{C}\right)   &  =1_{C}\ \ \ \ \ \ \ \ \ \ \text{and}\\
f\left(  1_{C}\right)   &  =1_{C}\ \ \ \ \ \ \ \ \ \ \text{and}\\
\operatorname*{Prim}C  &  \subseteq\operatorname*{Ker}\left(  e-f\right)
\ \ \ \ \ \ \ \ \ \ \text{and}\\
f\circ e  &  =e\circ f.
\end{align*}

Then:

\textbf{(a)} For any integer $u>1$, we have%
\[
\left(  e-f\right)  ^{u-1}\left(  C_{\leq u}\right)  \subseteq
\operatorname*{Prim}C.
\]

\textbf{(b)} For any positive integer $u$, we have%
\[
\left(  e-f\right)  ^{u}\left(  C_{\leq u}\right)  =0.
\]

\end{corollary}

See Section \ref{sec.pf.cor.id-f.cfc} for a proof of this corollary.

The particular case of Corollary \ref{cor.id-f.cfc1} for $e=\operatorname*{id}%
$ is particularly simple:

\begin{corollary}
\label{cor.id-f.cfc1id}Let $C$ be a connected filtered $\mathbf{k}$-coalgebra
with filtration $\left(  C_{\leq0},C_{\leq1},C_{\leq2},\ldots\right)  $.

Let $f:C\rightarrow C$ be a $\mathbf{k}$-coalgebra homomorphism such that
\[
f\left(  1_{C}\right)  =1_{C}\ \ \ \ \ \ \ \ \ \ \text{and}%
\ \ \ \ \ \ \ \ \ \ \operatorname*{Prim}C\subseteq\operatorname*{Ker}\left(
\operatorname*{id}-f\right)  .
\]

Then:

\textbf{(a)} For any integer $u>1$, we have%
\[
\left(  \operatorname*{id}-f\right)  ^{u-1}\left(  C_{\leq u}\right)
\subseteq\operatorname*{Prim}C.
\]

\textbf{(b)} For any positive integer $u$, we have%
\[
\left(  \operatorname*{id}-f\right)  ^{u}\left(  C_{\leq u}\right)  =0.
\]

\end{corollary}

Again, the proof of this corollary can be found in Section
\ref{sec.pf.cor.id-f.cfc}.

Note that Corollary \ref{cor.id-f.cfc1id} \textbf{(b)} is precisely
\cite[Theorem 37.1 \textbf{(a)}]{logid}.

\subsection{\label{sec.filhopf}Connected filtered bialgebras and Hopf
algebras}

We shall now apply our above results to connected filtered bialgebras and Hopf
algebras. We first define what we mean by these notions:

\begin{definition}
\label{def.fil-bial}\textbf{(a)} A \emph{filtered }$\mathbf{k}$%
\emph{-bialgebra} means a $\mathbf{k}$-bialgebra $H$ equipped with an infinite
sequence $\left(  H_{\leq0},H_{\leq1},H_{\leq2},\ldots\right)  $ of
$\mathbf{k}$-submodules of $H$ satisfying the following five conditions:

\begin{itemize}
\item We have%
\[
H_{\leq0}\subseteq H_{\leq1}\subseteq H_{\leq2}\subseteq\cdots.
\]

\item We have%
\[
\bigcup_{n\in\mathbb{N}}H_{\leq n}=H.
\]

\item We have
\[
\Delta\left(  H_{\leq n}\right)  \subseteq\sum_{i=0}^{n}H_{\leq i}\otimes
H_{\leq n-i}\ \ \ \ \ \ \ \ \ \ \text{for each }n\in\mathbb{N}.
\]
(Here, the \textquotedblleft$H_{\leq i}\otimes H_{\leq n-i}$\textquotedblright%
\ on the right hand side means the image of $H_{\leq i}\otimes H_{\leq n-i}$
under the canonical map $H_{\leq i}\otimes H_{\leq n-i}\rightarrow H\otimes H$
that is obtained by tensoring the two inclusion maps $H_{\leq i}\rightarrow H$
and $H_{\leq n-i}\rightarrow H$ together.)

\item We have $H_{\leq i}H_{\leq j}\subseteq H_{\leq i+j}$ for any
$i,j\in\mathbb{N}$. (Here, $H_{\leq i}H_{\leq j}$ denotes the $\mathbf{k}%
$-linear span of the set of all products $ab$ with $a\in H_{\leq i}$ and $b\in
H_{\leq j}$.)

\item The unity of the $\mathbf{k}$-algebra $H$ belongs to $H_{\leq0}$.{}
\end{itemize}

The sequence $\left(  H_{\leq0},H_{\leq1},H_{\leq2},\ldots\right)  $ is called
the \emph{filtration} of the filtered $\mathbf{k}$-bialgebra $H$.

\textbf{(b)} A \emph{filtered }$\mathbf{k}$\emph{-Hopf algebra} means a
filtered $\mathbf{k}$-bialgebra $H$ such that the $\mathbf{k}$-bialgebra $H$
is a Hopf algebra (i.e., has an antipode) and such that the antipode $S$ of
$H$ respects the filtration (i.e., satisfies $S\left(  H_{\leq n}\right)
\subseteq H_{\leq n}$ for each $n\in\mathbb{N}$).
\end{definition}

The $H_{\leq i}H_{\leq j}\subseteq H_{\leq i+j}$ condition in Definition
\ref{def.fil-bial} \textbf{(a)} will not actually be used in what follows.
Thus, we could have omitted it; but this would have resulted in a less common
(and less well-behaved in other ways) concept of \textquotedblleft filtered
bialgebra\textquotedblright. Likewise, we have included the $S\left(  H_{\leq
n}\right)  \subseteq H_{\leq n}$ condition in Definition \ref{def.fil-bial}
\textbf{(b)}, even though we will never use it.

Every $\mathbf{k}$-bialgebra is automatically a $\mathbf{k}$-coalgebra. Thus,
every filtered $\mathbf{k}$-bialgebra is automatically a filtered $\mathbf{k}%
$-coalgebra. This allows the following definition:

\begin{definition}
A filtered $\mathbf{k}$-bialgebra $H$ is said to be \emph{connected} if the
filtered $\mathbf{k}$-coalgebra $H$ is connected.
\end{definition}

Thus, if $H$ is a connected filtered $\mathbf{k}$-bialgebra, then Definition
\ref{def.con-fil-coal} \textbf{(b)} defines a \textquotedblleft
unity\textquotedblright\ $1_{H}$ of $H$. This appears to cause an awkward
notational quandary, since $H$ already has a unity by virtue of being a
$\mathbf{k}$-algebra (and this latter unity is also commonly denoted by
$1_{H}$). Fortunately, this cannot cause any confusion, since these two
unities are identical, as the following proposition shows:

\begin{proposition}
\label{prop.fil-bial.1=1}Let $H$ be a connected filtered $\mathbf{k}%
$-bialgebra. Then, the unity $1_{H}$ defined according to Definition
\ref{def.con-fil-coal} \textbf{(b)} equals the unity of the $\mathbf{k}%
$-algebra $H$.
\end{proposition}

\begin{proof}
[Proof of Proposition \ref{prop.fil-bial.1=1}.]Both unities in question belong
to $H_{\leq0}$ (indeed, the former does so by its definition, whereas the
latter does so because $H$ is a filtered $\mathbf{k}$-bialgebra) and are sent
to $1_{\mathbf{k}}$ by the map $\epsilon$ (indeed, the former does so by its
definition, whereas the latter does so by the axioms of a $\mathbf{k}%
$-bialgebra). However, since the map $\epsilon\mid_{H_{\leq0}}$ is a
$\mathbf{k}$-module isomorphism (because the filtered $\mathbf{k}$-coalgebra
$H$ is connected), these two properties uniquely determine these unities.
Thus, these two unities are equal. Proposition \ref{prop.fil-bial.1=1} is thus proven.
\end{proof}

In Definition \ref{def.con-fil-coal}, we have defined the notion of a
\textquotedblleft primitive element\textquotedblright\ of a connected filtered
$\mathbf{k}$-coalgebra $C$. In the same way, we can define a \textquotedblleft
primitive element\textquotedblright\ of a $\mathbf{k}$-bialgebra $H$ (using
the unity of the $\mathbf{k}$-algebra $H$ instead of $1_{C}$):

\begin{definition}
\label{def.prim-in-bialg}Let $H$ be a $\mathbf{k}$-bialgebra with unity
$1_{H}$.

\textbf{(a)} An element $x$ of $H$ is said to be \emph{primitive} if
$\Delta\left(  x\right)  =x\otimes1_{H}+1_{H}\otimes x$.

\textbf{(b)} The set of all primitive elements of $H$ is denoted by
$\operatorname*{Prim}H$.
\end{definition}

When $H$ is a connected filtered $\mathbf{k}$-bialgebra, Definition
\ref{def.prim-in-bialg} \textbf{(a)}\ agrees with Definition
\ref{def.con-fil-coal} \textbf{(c)}, since Proposition \ref{prop.fil-bial.1=1}
shows that the two meanings of $1_{H}$ are actually identical. Thus, when $H$
is a connected filtered $\mathbf{k}$-bialgebra, Definition
\ref{def.prim-in-bialg} \textbf{(b)}\ agrees with Definition
\ref{def.con-fil-coal} \textbf{(d)}. The notation $\operatorname*{Prim}H$ is
therefore unambiguous.

For later use, we state some basic properties of the antipode in a Hopf algebra:

\begin{lemma}
\label{lem.bialg.antip-props}Let $H$ be a $\mathbf{k}$-Hopf algebra with unity
$1_{H}\in H$ and antipode $S\in\operatorname*{End}H$. Then:

\textbf{(a)} The map $S^{2}:H\rightarrow H$ is a $\mathbf{k}$-coalgebra homomorphism.

\textbf{(b)} We have $S\left(  1_{H}\right)  =1_{H}$.

\textbf{(c)} We have $S\left(  x\right)  =-x$ for every primitive element $x$
of $H$.

\textbf{(d)} We have $S^{2}\left(  x\right)  =x$ for every primitive element
$x$ of $H$.
\end{lemma}

All of these facts are easy to prove and well-known (see \cite[proof of Lemma
37.8]{logid} for detailed proofs, or derive them easily from \cite[Proposition
1.2.17 (1)]{HecSch20}).

We can now state our main consequence for connected filtered Hopf algebras;
all proofs can be found in Section \ref{sec.pf.filhopf} below:

\begin{corollary}
\label{cor.id-S2.p}Let $H$ be a connected filtered $\mathbf{k}$-Hopf algebra
with antipode $S$ and
filtration $\left(  H_{\leq0},H_{\leq1},H_{\leq2},\ldots\right)$.

Let $p$ be a positive integer such that%
\begin{equation}
\left(  \operatorname*{id}-S^{2}\right)  \left(  H_{\leq p}\right)  =0.
\label{eq.cor.id-S2.p.ass-ann}%
\end{equation}

Then:

\textbf{(a)} For any integer $u>p$, we have%
\begin{equation}
\left(  \operatorname*{id}-S^{2}\right)  ^{u-p}\left(  H_{\leq u}\right)
\subseteq\operatorname*{Prim}H \label{eq.cor.id-S2.p.claim1}%
\end{equation}
and%
\begin{equation}
\left(  \left(  \operatorname*{id}+S\right)  \circ\left(  \operatorname*{id}%
-S^{2}\right)  ^{u-p}\right)  \left(  H_{\leq u}\right)  =0.
\label{eq.cor.id-S2.p.claim3}%
\end{equation}

\textbf{(b)} For any integer $u\geq p$, we have%
\begin{equation}
\left(  \operatorname*{id}-S^{2}\right)  ^{u-p+1}\left(  H_{\leq u}\right)
=0. \label{eq.cor.id-S2.p.claim2}%
\end{equation}

\end{corollary}

Specializing this to $p=1$, we can easily obtain the following:

\begin{corollary}
\label{cor.id-S2.1}Let $H$ be a connected filtered $\mathbf{k}$-Hopf algebra
with antipode $S$ and filtration $\left(  H_{\leq0},H_{\leq1},H_{\leq2},\ldots\right)  $.
Then:

\textbf{(a)} For any integer $u>1$, we have%
\begin{equation}
\left(  \operatorname*{id}-S^{2}\right)  ^{u-1}\left(  H_{\leq u}\right)
\subseteq\operatorname*{Prim}H \label{eq.cor.id-S2.1.claim1}%
\end{equation}
and%
\begin{equation}
\left(  \left(  \operatorname*{id}+S\right)  \circ\left(  \operatorname*{id}%
-S^{2}\right)  ^{u-1}\right)  \left(  H_{\leq u}\right)  =0.
\label{eq.cor.id-S2.1.claim3}%
\end{equation}

\textbf{(b)} For any positive integer $u$, we have%
\begin{equation}
\left(  \operatorname*{id}-S^{2}\right)  ^{u}\left(  H_{\leq u}\right)  =0.
\label{eq.cor.id-S2.1.claim2}%
\end{equation}

\end{corollary}

Corollary \ref{cor.id-S2.1} \textbf{(b)} has already appeared in \cite[Theorem
37.7 \textbf{(a)}]{logid}.

\subsection{\label{sec.hopf}Connected graded Hopf algebras}

Let us now specialize our results even further to connected \textbf{graded}
Hopf algebras. We have already seen that any graded $\mathbf{k}$-coalgebra
automatically becomes a filtered $\mathbf{k}$-coalgebra. In the same way, any
graded $\mathbf{k}$-Hopf algebra automatically becomes a filtered $\mathbf{k}%
$-Hopf algebra. Moreover, a graded $\mathbf{k}$-Hopf algebra $H$ is connected
(in the sense that $H_{0}\cong\mathbf{k}$ as $\mathbf{k}$-modules) if and only
if the filtered $\mathbf{k}$-coalgebra $H$ is connected. (This follows easily
from \cite[Exercise 1.3.20 (e)]{GriRei}.) Thus, our above results for
connected filtered $\mathbf{k}$-Hopf algebras can be applied to connected
graded $\mathbf{k}$-Hopf algebras. From Corollary \ref{cor.id-S2.1}, we easily
obtain the following:

\begin{corollary}
\label{cor.id-S2.gr1}Let $H$ be a connected graded $\mathbf{k}$-Hopf algebra
with antipode $S$. Then, for any positive integer $u$, we have%
\begin{equation}
\left(  \operatorname*{id}-S^{2}\right)  ^{u-1}\left(  H_{u}\right)
\subseteq\operatorname*{Prim}H \label{eq.cor.id-S2.gr1.claim1}%
\end{equation}
and%
\begin{equation}
\left(  \left(  \operatorname*{id}+S\right)  \circ\left(  \operatorname*{id}%
-S^{2}\right)  ^{u-1}\right)  \left(  H_{u}\right)  =0
\label{eq.cor.id-S2.gr1.claim3}%
\end{equation}
and%
\begin{equation}
\left(  \operatorname*{id}-S^{2}\right)  ^{u}\left(  H_{u}\right)  =0.
\label{eq.cor.id-S2.gr1.claim2}%
\end{equation}

\end{corollary}

We will prove this corollary -- as well as all others stated in Section
\ref{sec.hopf} -- in Section \ref{sec.pf.hopf} further below. In essence,
Corollary \ref{cor.id-S2.gr1} follows from Corollary \ref{cor.id-S2.1}, once
the case $u=1$ (which is not covered by Corollary \ref{cor.id-S2.1}
\textbf{(a)}) is accounted for.

The equality (\ref{eq.cor.id-S2.gr1.claim3}) in Corollary \ref{cor.id-S2.gr1}
yields \cite[Lemma 12.50]{Aguiar17}, whereas the equality
(\ref{eq.cor.id-S2.gr1.claim2}) yields \cite[Proposition 7]{AguLau14}. Next,
we apply Corollary \ref{cor.id-S2.p} to the graded setting:

\begin{corollary}
\label{cor.id-S2.grp}Let $H$ be a connected graded $\mathbf{k}$-Hopf algebra
with antipode $S$.

Let $p$ be a positive integer such that all $i\in\left\{  2,3,\ldots
,p\right\}  $ satisfy%
\begin{equation}
\left(  \operatorname*{id}-S^{2}\right)  \left(  H_{i}\right)  =0.
\label{eq.cor.id-S2.grp.ass-ann}%
\end{equation}

Then:

\textbf{(a)} For any integer $u>p$, we have%
\begin{equation}
\left(  \operatorname*{id}-S^{2}\right)  ^{u-p}\left(  H_{\leq u}\right)
\subseteq\operatorname*{Prim}H \label{eq.cor.id-S2.grp.claim1}%
\end{equation}
and%
\begin{equation}
\left(  \left(  \operatorname*{id}+S\right)  \circ\left(  \operatorname*{id}%
-S^{2}\right)  ^{u-p}\right)  \left(  H_{\leq u}\right)  =0.
\label{eq.cor.id-S2.grp.claim3}%
\end{equation}

\textbf{(b)} For any integer $u\geq p$, we have%
\begin{equation}
\left(  \operatorname*{id}-S^{2}\right)  ^{u-p+1}\left(  H_{\leq u}\right)
=0. \label{eq.cor.id-S2.grp.claim2}%
\end{equation}

\end{corollary}

The particular case of Corollary \ref{cor.id-S2.grp} for $p=2$ is the most
useful, as the condition (\ref{eq.cor.id-S2.grp.ass-ann}) boils down to the
equality $\left(  \operatorname*{id}-S^{2}\right)  \left(  H_{2}\right)  =0$
in this case, and the latter equality is satisfied rather frequently. Here is
one sufficient criterion:

\begin{corollary}
\label{cor.id-S2.gr2}Let $H$ be a connected graded $\mathbf{k}$-Hopf algebra
with antipode $S$. Assume that%
\begin{equation}
ab=ba\ \ \ \ \ \ \ \ \ \ \text{for every }a,b\in H_{1}.
\label{eq.cor.id-S2.gr2.commH1}%
\end{equation}
Then:

\textbf{(a)} We have
\[
\left(  \operatorname*{id}-S^{2}\right)  \left(  H_{2}\right)  =0.
\]

\textbf{(b)} For any integer $u>2$, we have%
\begin{equation}
\left(  \operatorname*{id}-S^{2}\right)  ^{u-2}\left(  H_{\leq u}\right)
\subseteq\operatorname*{Prim}H \label{eq.cor.id-S2.gr2.b.claim1}%
\end{equation}
and%
\begin{equation}
\left(  \left(  \operatorname*{id}+S\right)  \circ\left(  \operatorname*{id}%
-S^{2}\right)  ^{u-2}\right)  \left(  H_{\leq u}\right)  =0.
\label{eq.cor.id-S2.gr2.b.claim3}%
\end{equation}

\textbf{(c)} For any integer $u>1$, we have%
\begin{equation}
\left(  \operatorname*{id}-S^{2}\right)  ^{u-1}\left(  H_{\leq u}\right)  =0.
\label{eq.cor.id-S2.gr2.b.claim2}%
\end{equation}

\end{corollary}

The equality (\ref{eq.cor.id-S2.gr2.b.claim2}) generalizes \cite[Example
8]{AguLau14}. Indeed, if $H$ is the Malvenuto--Reutenauer Hopf
algebra\footnote{See \cite[\S 12.1]{Meliot17}, \cite[\S 7.1]{HaGuKi10} or
\cite[\S 8.1]{GriRei} for the definition of this Hopf algebra. (It is denoted
$\operatorname*{FQSym}$ in \cite{Meliot17} and \cite{GriRei}, and denoted
$MPR$ in \cite{HaGuKi10}.)}, then the condition (\ref{eq.cor.id-S2.gr2.commH1}%
) is satisfied (since $H_{1}$ is a free $\mathbf{k}$-module of rank $1$ in
this case); therefore, Corollary \ref{cor.id-S2.gr2} \textbf{(c)} can be
applied in this case, and we recover \cite[Example 8]{AguLau14}. Likewise, we
can obtain the same result if $H$ is the Hopf algebra $\operatorname*{WQSym}$
of word quasisymmetric functions\footnote{See (e.g.) \cite[\S 4.3.2]{MeNoTh13}
for a definition of this Hopf algebra.}.

It is worth noticing that the condition (\ref{eq.cor.id-S2.gr2.commH1}) is
only sufficient, but not necessary for (\ref{eq.cor.id-S2.gr2.b.claim2}). For
example, if $H$ is the tensor algebra of a free $\mathbf{k}$-module $V$ of
rank $\geq2$, then (\ref{eq.cor.id-S2.gr2.b.claim2}) holds (since $H$ is
cocommutative, so that $S^{2}=\operatorname*{id}$), but
(\ref{eq.cor.id-S2.gr2.commH1}) does not (since $u\otimes v\neq v\otimes u$ if
$u$ and $v$ are two distinct basis vectors of $V$).

An example of a connected graded Hopf algebra $H$ that does \textbf{not}
satisfy (\ref{eq.cor.id-S2.gr2.b.claim2}) (and thus does not satisfy
(\ref{eq.cor.id-S2.gr2.commH1}) either) is not hard to construct:

\begin{example}
Assume that the ring $\mathbf{k}$ is not trivial. Let $H$ be the free
$\mathbf{k}$-algebra with three generators $a,b,c$. We equip this $\mathbf{k}%
$-algebra $H$ with a grading, by requiring that its generators $a,b,c$ are
homogeneous of degrees $1,1,2$, respectively. Next, we define a
comultiplication $\Delta$ on $H$ by setting
\begin{align*}
\Delta\left(  a\right)   &  =a\otimes1+1\otimes a;\\
\Delta\left(  b\right)   &  =b\otimes1+1\otimes b;\\
\Delta\left(  c\right)   &  =c\otimes1+a\otimes b+1\otimes c
\end{align*}
(where $1$ is the unity of $H$). Furthermore, we define a counit $\epsilon$ on
$H$ by setting
\[
    \epsilon\left(  a\right)  =\epsilon\left(  b\right) =\epsilon\left(  c\right)  =0.
\]
It is straightforward to see that $H$ thus
becomes a connected graded $\mathbf{k}$-bialgebra, hence (by \cite[Proposition
1.4.16]{GriRei}) a connected graded $\mathbf{k}$-Hopf algebra. Its antipode
$S$ is easily seen to satisfy $S\left(  c\right)  =ab-c$ and $S^{2}\left(
c\right)  =ba-ab+c\neq c$; thus, $\left(  \operatorname*{id}-S^{2}\right)
\left(  H_{2}\right)  \neq0$. Hence, (\ref{eq.cor.id-S2.gr2.b.claim2}) does
not hold for $u=2$.
\end{example}

The Hopf algebra $H$ in this example is in fact an instance of a general
construction (to be elaborated upon in future work) of connected graded
$\mathbf{k}$-Hopf algebras that are \textquotedblleft
generic\textquotedblright\ (in the sense that their structure maps satisfy no
relations other than ones that hold in every connected graded $\mathbf{k}%
$-Hopf algebra).

\begin{remark}
A brave reader might wonder whether the connectedness condition in Corollary
\ref{cor.id-S2.gr1} could be replaced by something weaker -- e.g., instead of
requiring $H$ to be connected, we might require that the subalgebra $H_{0}$ be
commutative. However, such a requirement would be insufficient. In fact, let
$\mathbf{k}=\mathbb{C}$. Then, for any integer $n>1$ and any primitive $n$-th
root of unity $q\in\mathbf{k}$, the Taft algebra $H_{n,q}$ defined in
\cite[\S 7.3]{Radfor12} can be viewed as a graded Hopf algebra (with $a\in
H_{0}$ and $x\in H_{1}$) whose subalgebra $H_{0}=\mathbf{k}\left[  a\right]
/\left(  a^{n}-1\right)  $ is commutative, but whose antipode $S$ does not
satisfy $\left(  \operatorname*{id}-S^{2}\right)  ^{k}\left(  H_{1}\right)
=0$ for any $k\in\mathbb{N}$ (since $S^{2}\left(  x\right)  =q^{-1}x$ and
therefore $\left(  \operatorname*{id}-S^{2}\right)  ^{k}\left(  x\right)
=\left(  1-q^{-1}\right)  ^{k}x\neq0$ because $q^{-1}\neq1$).
\end{remark}

\section{Proofs}

We shall now prove all statements left unproved above.

\subsection{\label{sec.pf.thm.id-f.gen}Proof of Theorem \ref{thm.id-f.gen}}

\begin{proof}
[Proof of Theorem \ref{thm.id-f.gen}.]We shall prove
(\ref{eq.thm.id-f.gen.clm1}) and (\ref{eq.thm.id-f.gen.clm2}) by strong
induction on $u$:

\textit{Induction step:} Fix an integer $n>p$. Assume (as the induction
hypothesis) that (\ref{eq.thm.id-f.gen.clm1}) and (\ref{eq.thm.id-f.gen.clm2})
hold for all integers $u>p$ satisfying $u<n$. We must prove that
(\ref{eq.thm.id-f.gen.clm1}) and (\ref{eq.thm.id-f.gen.clm2}) hold for $u=n$.
In other words, we must prove that%
\[
\left(  e-f\right)  ^{n-p}\left(  D_{n}\right)  \subseteq\operatorname*{Ker}%
\delta
\]
and%
\[
\left(  e-f\right)  ^{n-p+1}\left(  D_{n}\right)  =0.
\]
We shall focus on proving the first of these two equalities; the second will
then easily follow from (\ref{pf.thm.id-f.gen.ass-Ker}).

Consider the $\mathbf{k}$-algebras $\operatorname*{End}D$ and
$\operatorname*{End}\left(  D\otimes D\right)  $. (The multiplication in each
of these $\mathbf{k}$-algebras is composition of $\mathbf{k}$-linear maps.)
Note that $u\otimes v\in\operatorname*{End}\left(  D\otimes D\right)  $ for
any $u,v\in\operatorname*{End}D$.

We have $e,f\in\operatorname*{End}D$. Let us define two elements
$g\in\operatorname*{End}D$ and $h\in\operatorname*{End}\left(  D\otimes
D\right)  $ by%
\[
g=e-f\ \ \ \ \ \ \ \ \ \ \text{and}\ \ \ \ \ \ \ \ \ \ h=e\otimes e-f\otimes
f.
\]
Then, we easily obtain%
\[
h=g\otimes f+e\otimes g.
\]
Moreover, (\ref{pf.thm.id-f.gen.ass-ann}) rewrites as $g\left(  D_{1}%
+D_{2}+\cdots+D_{p}\right)  =0$ (since $g=e-f$). Thus,%
\begin{equation}
g\left(  D_{u}\right)  =0\ \ \ \ \ \ \ \ \ \ \text{for all }u\in\left\{
1,2,\ldots,p\right\}  . \label{pf.thm.id-f.gen.gDu=0}%
\end{equation}

Now, recall that the multiplication in the $\mathbf{k}$-algebra
$\operatorname*{End}D$ is composition of maps. Thus, $\alpha\beta=\alpha
\circ\beta$ for any $\alpha,\beta\in\operatorname*{End}D$. (The same holds for
$\operatorname*{End}\left(  D\otimes D\right)  $.) Hence,
(\ref{pf.thm.id-f.gen.ass-comm}) rewrites as $fe=ef$. In other words, $e$ and
$f$ commute. Therefore, the difference $g=e-f$ commutes with $e$ and $f$ as
well, i.e., we have $ge=eg$ and $gf=fg$. Hence, for each $i\in\mathbb{N}$ and
$j\in\mathbb{N}$, we have%
\begin{equation}
g^{i}e^{j}=e^{j}g^{i} \label{pf.thm.id-f.gen.giej}%
\end{equation}
(since powers of commuting elements always commute).

Furthermore, in $\operatorname*{End}\left(  D\otimes D\right)  $, we have%
\[
\left(  g\otimes f\right)  \left(  e\otimes g\right)  =\underbrace{\left(
ge\right)  }_{=eg}\otimes\underbrace{\left(  fg\right)  }_{=gf}=\left(
eg\right)  \otimes\left(  gf\right)  =\left(  e\otimes g\right)  \left(
g\otimes f\right)  .
\]
Hence, we can apply the binomial formula to $g\otimes f$ and $e\otimes g$. We
thus conclude that each $k\in\mathbb{N}$ satisfies%
\begin{align}
\left(  g\otimes f+e\otimes g\right)  ^{k}  &  =\sum_{r=0}^{k}\dbinom{k}%
{r}\underbrace{\left(  g\otimes f\right)  ^{r}\left(  e\otimes g\right)
^{k-r}}_{\substack{=\left(  g^{r}\otimes f^{r}\right)  \left(  e^{k-r}\otimes
g^{k-r}\right)  \\=\left(  g^{r}e^{k-r}\right)  \otimes\left(  f^{r}%
g^{k-r}\right)  }}=\sum_{r=0}^{k}\dbinom{k}{r}\underbrace{\left(  g^{r}%
e^{k-r}\right)  }_{\substack{=e^{k-r}g^{r}\\\text{(by
(\ref{pf.thm.id-f.gen.giej}))}}}\otimes\left(  f^{r}g^{k-r}\right) \nonumber\\
&  =\sum_{r=0}^{k}\dbinom{k}{r}\underbrace{\left(  e^{k-r}g^{r}\right)
\otimes\left(  f^{r}g^{k-r}\right)  }_{=\left(  e^{k-r}\otimes f^{r}\right)
\left(  g^{r}\otimes g^{k-r}\right)  }\nonumber\\
&  =\sum_{r=0}^{k}\dbinom{k}{r}\left(  e^{k-r}\otimes f^{r}\right)  \left(
g^{r}\otimes g^{k-r}\right)  . \label{pf.thm.id-f.gen.ass-Ker.binom}%
\end{align}

For each $k\in\mathbb{N}$ and $r\in\mathbb{N}$, we define a map $h_{k,r}%
\in\operatorname*{End}\left(  D\otimes D\right)  $ by%
\begin{equation}
h_{k,r}=\dbinom{k}{r}\left(  e^{k-r}\otimes f^{r}\right)  \left(  g^{r}\otimes
g^{k-r}\right)  . \label{pf.thm.id-f.gen.hkr=}%
\end{equation}
Thus, we can rewrite (\ref{pf.thm.id-f.gen.ass-Ker.binom}) as follows: Each
$k\in\mathbb{N}$ satisfies%
\begin{equation}
h^{k}=\sum_{r=0}^{k}h_{k,r} \label{pf.thm.id-f.gen.hk=}%
\end{equation}
(because of (\ref{pf.thm.id-f.gen.hkr=}), and because $h=g\otimes f+e\otimes
g$).

Subtracting (\ref{pf.thm.id-f.gen.ass-mor}) from
(\ref{pf.thm.id-f.gen.ass-mor'}), we obtain
\[
\left(  e\otimes e-f\otimes f\right)  \circ\delta=\delta\circ\left(
e-f\right)
\]
(since $\delta$ is $\mathbf{k}$-linear). In view of $g=e-f$ and $h=e\otimes
e-f\otimes f$, this rewrites as $h\circ\delta=\delta\circ g$. In other words,
$\delta\circ g=h\circ\delta$. Hence, by induction on $k$, we easily see that%
\begin{equation}
\delta\circ g^{k}=h^{k}\circ\delta\ \ \ \ \ \ \ \ \ \ \text{for each }%
k\in\mathbb{N}. \label{pf.thm.id-f.gen.ass-Ker.delgk=}%
\end{equation}

Our induction hypothesis says that (\ref{eq.thm.id-f.gen.clm1}) and
(\ref{eq.thm.id-f.gen.clm2}) hold for all integers $u>p$ satisfying $u<n$. In
particular, (\ref{eq.thm.id-f.gen.clm2}) holds for all integers $u>p$
satisfying $u<n$. In other words, for each integer $u>p$ satisfying $u<n$, we
have
\begin{equation}
g^{u-p+1}\left(  D_{u}\right)  =0 \label{pf.thm.id-f.gen.ass-Ker.IH3}%
\end{equation}
(since $g=e-f$). Hence, it is easy to see that every positive integer $u<n$
and every positive integer $v>u-p$ satisfy%
\begin{equation}
g^{v}\left(  D_{u}\right)  =0. \label{pf.thm.id-f.gen.ass-Ker.IH4}%
\end{equation}
(Indeed, if $u>p$, then this follows from (\ref{pf.thm.id-f.gen.ass-Ker.IH3}),
because $v\geq u-p+1$. However, if $u\leq p$, then
(\ref{pf.thm.id-f.gen.ass-Ker.IH4}) follows from (\ref{pf.thm.id-f.gen.gDu=0}%
), because $v\geq1$. Thus, (\ref{pf.thm.id-f.gen.ass-Ker.IH4}) is proved in
all possible cases.)

Now, let $k=n-p$. Then, $k>0$ (since $n>p$), so that $k\in\mathbb{N}$.
Furthermore, (\ref{pf.thm.id-f.gen.ass-Ker.delgk=}) yields $\delta\circ
g^{k}=h^{k}\circ\delta$. Thus,
\begin{align}
\left(  \delta\circ g^{k}\right)  \left(  D_{n}\right)   &  =\left(
h^{k}\circ\delta\right)  \left(  D_{n}\right)  =h^{k}\left(  \delta\left(
D_{n}\right)  \right) \nonumber\\
&  \subseteq h^{k}\left(  \sum_{i=1}^{n-1}D_{i}\otimes D_{n-i}\right)
\ \ \ \ \ \ \ \ \ \ \left(  \text{by (\ref{pf.thm.id-f.gen.ass-gr})}\right)
\nonumber\\
&  =\sum_{i=1}^{n-1}h^{k}\left(  D_{i}\otimes D_{n-i}\right)  .
\label{pf.thm.id-f.gen.in-bigsum}%
\end{align}

We shall now prove that each $i\in\left\{  1,2,\ldots,n-1\right\}  $ and each
$r\in\left\{  0,1,\ldots,k\right\}  $ satisfy
\begin{equation}
\left(  g^{r}\otimes g^{k-r}\right)  \left(  D_{i}\otimes D_{n-i}\right)  =0.
\label{pf.thm.id-f.gen.ass-Ker.tens=0}%
\end{equation}

[\textit{Proof of (\ref{pf.thm.id-f.gen.ass-Ker.tens=0}):} Fix $i\in\left\{
1,2,\ldots,n-1\right\}  $ and $r\in\left\{  0,1,\ldots,k\right\}  $. We must
prove (\ref{pf.thm.id-f.gen.ass-Ker.tens=0}).

We have $i\in\left\{  1,2,\ldots,n-1\right\}  $. Thus, $0<i<n$, so that
$0<n-i<n$. Also, $\min\left\{  k,i\right\}  >0$ (since $k>0$ and $i>0$). Also,
$k=n-p>i-p$ (since $n>i$) and $i>i-p$ (since $p>0$). Therefore, $\min\left\{
k,i\right\}  >i-p$.

We are in one of the following two cases:

\textit{Case 1:} We have $r\geq\min\left\{  k,i\right\}  $.

\textit{Case 2:} We have $r<\min\left\{  k,i\right\}  $.

Let us first consider Case 1. In this case, we have $r\geq\min\left\{
k,i\right\}  $. This entails that $r\geq\min\left\{  k,i\right\}  >i-p$. Moreover,
the integer $r$ is positive (since $r\geq\min\left\{  k,i\right\}  >0$).
Hence, (\ref{pf.thm.id-f.gen.ass-Ker.IH4}) (applied to $u=i$ and $v=r$) yields
$g^{r}\left(  D_{i}\right)  =0$ (since $r>i-p$). Now,%
\[
\left(  g^{r}\otimes g^{k-r}\right)  \left(  D_{i}\otimes D_{n-i}\right)
=\underbrace{g^{r}\left(  D_{i}\right)  }_{=0}\otimes g^{k-r}\left(
D_{n-i}\right)  =0.
\]
Thus, (\ref{pf.thm.id-f.gen.ass-Ker.tens=0}) is proved in Case 1.

Next, let us consider Case 2. In this case, we have $r<\min\left\{
k,i\right\}  $. Now, the integer $k-r$ is positive (since $r<\min\left\{
k,i\right\}  \leq k$). Furthermore, from $k=n-p$, we obtain%
\[
k-r=n-p-\underbrace{r}_{<\min\left\{  k,i\right\}  \leq i}>n-p-i=n-i-p.
\]
Hence, (\ref{pf.thm.id-f.gen.ass-Ker.IH4}) (applied to $u=n-i$ and $v=k-r$)
yields $g^{k-r}\left(  D_{n-i}\right)  =0$. Now,%
\[
\left(  g^{r}\otimes g^{k-r}\right)  \left(  D_{i}\otimes D_{n-i}\right)
=g^{r}\left(  D_{i}\right)  \otimes\underbrace{g^{k-r}\left(  D_{n-i}\right)
}_{=0}=0.
\]
Thus, (\ref{pf.thm.id-f.gen.ass-Ker.tens=0}) is proved in Case 2.

The proof of (\ref{pf.thm.id-f.gen.ass-Ker.tens=0}) is now complete.]

Using (\ref{pf.thm.id-f.gen.ass-Ker.tens=0}), we can easily see that each
$i\in\left\{  1,2,\ldots,n-1\right\}  $ and each $r\in\left\{  0,1,\ldots
,k\right\}  $ satisfy
\begin{equation}
h_{k,r}\left(  D_{i}\otimes D_{n-i}\right)  =0.
\label{pf.thm.id-f.gen.ass-Ker.tens=02}%
\end{equation}
(Indeed, (\ref{pf.thm.id-f.gen.hkr=}) expresses $h_{k,r}$ as a product, but
the rightmost factor $g^{r}\otimes g^{k-r}$ in this product already kills
$D_{i}\otimes D_{n-i}$ according to (\ref{pf.thm.id-f.gen.ass-Ker.tens=0}).)

Now, each $i\in\left\{  1,2,\ldots,n-1\right\}  $ satisfies%
\begin{align}
h^{k}\left(  D_{i}\otimes D_{n-i}\right)   &  =\left(  \sum_{r=0}^{k}%
h_{k,r}\right)  \left(  D_{i}\otimes D_{n-i}\right)
\ \ \ \ \ \ \ \ \ \ \left(  \text{by (\ref{pf.thm.id-f.gen.hk=})}\right)
\nonumber\\
&  \subseteq\sum_{r=0}^{k}\underbrace{h_{k,r}\left(  D_{i}\otimes
D_{n-i}\right)  }_{\substack{=0\\\text{(by
(\ref{pf.thm.id-f.gen.ass-Ker.tens=02}))}}}=0.\nonumber
\end{align}
Hence, (\ref{pf.thm.id-f.gen.in-bigsum}) becomes%
\[
\left(  \delta\circ g^{k}\right)  \left(  D_{n}\right)  \subseteq\sum
_{i=1}^{n-1}\underbrace{h^{k}\left(  D_{i}\otimes D_{n-i}\right)  }%
_{\subseteq0}\subseteq0.
\]
In other words, $\delta\left(  g^{k}\left(  D_{n}\right)  \right)  \subseteq
0$. Equivalently,%
\[
g^{k}\left(  D_{n}\right)  \subseteq\operatorname*{Ker}\delta.
\]
Since $g=e-f$ and $k=n-p$, we can rewrite this as follows:%
\begin{equation}
\left(  e-f\right)  ^{n-p}\left(  D_{n}\right)  \subseteq\operatorname*{Ker}%
\delta. \label{pf.thm.id-f.gen.done1}%
\end{equation}

However, we can rewrite (\ref{pf.thm.id-f.gen.ass-Ker}) as $\left(
e-f\right)  \left(  \operatorname*{Ker}\delta\right)  =0$. Thus,%
\[
\left(  e-f\right)  ^{n-p+1}\left(  D_{n}\right)  =\left(  e-f\right)  \left(
\underbrace{\left(  e-f\right)  ^{n-p}\left(  D_{n}\right)  }%
_{\substack{\subseteq\operatorname*{Ker}\delta\\\text{(by
(\ref{pf.thm.id-f.gen.done1}))}}}\right)  \subseteq\left(  e-f\right)  \left(
\operatorname*{Ker}\delta\right)  =0.
\]
In other words,%
\begin{equation}
\left(  e-f\right)  ^{n-p+1}\left(  D_{n}\right)  =0.
\label{pf.thm.id-f.gen.done2}%
\end{equation}

We have now proved the relations (\ref{pf.thm.id-f.gen.done1}) and
(\ref{pf.thm.id-f.gen.done2}). In other words, (\ref{eq.thm.id-f.gen.clm1})
and (\ref{eq.thm.id-f.gen.clm2}) hold for $u=n$. This completes the induction
step. Thus, Theorem \ref{thm.id-f.gen} is proven.
\end{proof}

\subsection{\label{sec.pf.prop.cfc.delta2}Proof of Proposition
\ref{prop.cfc.delta2}}

Our next goal is to prove Proposition \ref{prop.cfc.delta2}. We shall work
towards this goal by proving a simple lemma:

\begin{lemma}
\label{lem.coalg.primitive-e0}Let $C$ be any $\mathbf{k}$-coalgebra. Let
$a,b,d\in C$ be three elements satisfying $\epsilon\left(  a\right)  =1$ and
$\epsilon\left(  b\right)  =1$ and $\Delta\left(  d\right)  =d\otimes
a+b\otimes d$. Then, $\epsilon\left(  d\right)  =0$.
\end{lemma}

We shall later apply Lemma \ref{lem.coalg.primitive-e0} to the case when
$a=b=1_{C}$ (and $C$ is either a connected filtered $\mathbf{k}$-coalgebra or
a $\mathbf{k}$-bialgebra, so that $1_{C}$ does make sense); however, it is not
any harder to prove it in full generality:

\begin{proof}
[Proof of Lemma \ref{lem.coalg.primitive-e0}.]Let $\gamma$ be the canonical
$\mathbf{k}$-module isomorphism $C\otimes\mathbf{k}\rightarrow C,\ c\otimes
\lambda\mapsto\lambda c$. One of the axioms of a coalgebra says that
$\gamma\circ\left(  \operatorname*{id}\otimes\epsilon\right)  \circ
\Delta=\operatorname*{id}$. Applying both sides of this equality to $d$, we
obtain $\left(  \gamma\circ\left(  \operatorname*{id}\otimes\epsilon\right)
\circ\Delta\right)  \left(  d\right)  =\operatorname*{id}\left(  d\right)
=d$. Hence,%
\begin{align*}
d  &  =\left(  \gamma\circ\left(  \operatorname*{id}\otimes\epsilon\right)
\circ\Delta\right)  \left(  d\right)  =\underbrace{\epsilon\left(  a\right)
}_{=1}d+\epsilon\left(  d\right)  b\ \ \ \ \ \ \ \ \ \ \left(  \text{since
}\Delta\left(  d\right)  =d\otimes a+b\otimes d\right) \\
&  =d+\epsilon\left(  d\right)  b.
\end{align*}
Subtracting $d$ from both sides, we obtain $\epsilon\left(  d\right)  b=0$.
Applying the map $\epsilon$ to this equality, we find $\epsilon\left(
\epsilon\left(  d\right)  b\right)  =0$. In view of%
\[
\epsilon\left(  \epsilon\left(  d\right)  b\right)  =\epsilon\left(  d\right)
\underbrace{\epsilon\left(  b\right)  }_{=1}=\epsilon\left(  d\right)  ,
\]
this rewrites as $\epsilon\left(  d\right)  =0$. This proves Lemma
\ref{lem.coalg.primitive-e0}.
\end{proof}

Next, let us define a \textquotedblleft reduced identity map\textquotedblright%
\ $\idbar$ for any connected filtered $\mathbf{k}$-coalgebra $C$, and explore
some of its properties:

\begin{lemma}
\label{lem.cfc.idbar}Let $C$ be a connected filtered $\mathbf{k}$-coalgebra
with filtration $\left(  C_{\leq0},C_{\leq1},C_{\leq2},\ldots\right)  $.
Define a $\mathbf{k}$-linear map $\idbar:C\rightarrow C$ by setting%
\[
\idbar\left(  c\right)  :=c-\epsilon\left(  c\right)  1_{C}%
\ \ \ \ \ \ \ \ \ \ \text{for each }c\in C.
\]
Define a $\mathbf{k}$-linear map $\delta:C\rightarrow C\otimes C$ by setting
\[
\delta\left(  c\right)  :=\Delta\left(  c\right)  -c\otimes1_{C}-1_{C}\otimes
c+\epsilon\left(  c\right)  1_{C}\otimes1_{C}\ \ \ \ \ \ \ \ \ \ \text{for
each }c\in C.
\]
Then:

\textbf{(a)} We have $\delta=\left(  \idbar\otimes\idbar\right)  \circ\Delta$.

\textbf{(b)} We have $\idbar\left(  C_{\leq n}\right)  \subseteq C_{\leq n}$
for each $n\in\mathbb{N}$.

\textbf{(c)} We have $\idbar\left(  C_{\leq0}\right)  =0$.
\end{lemma}

\begin{proof}
[Proof of Lemma \ref{lem.cfc.idbar}.]\textbf{(a)} Let $c\in C$. Write the
tensor $\Delta\left(  c\right)  \in C\otimes C$ in the form%
\begin{equation}
\Delta\left(  c\right)  =\sum_{i=1}^{m}c_{i}\otimes d_{i}
\label{pf.lem.cfc.idbar.a.1}%
\end{equation}
for some $m\in\mathbb{N}$, some $c_{1},c_{2},\ldots,c_{m}\in C$ and some
$d_{1},d_{2},\ldots,d_{m}\in C$.

According to the axioms of a coalgebra\footnote{Specifically, we are using the
axioms%
\[
\operatorname*{id}\nolimits_{C}=\gamma\circ\left(  \operatorname*{id}%
\otimes\epsilon\right)  \circ\Delta\ \ \ \ \ \ \ \ \ \ \text{and}%
\ \ \ \ \ \ \ \ \ \ \operatorname*{id}\nolimits_{C}=\gamma^{\prime}%
\circ\left(  \epsilon\otimes\operatorname*{id}\right)  \circ\Delta,
\]
where $\gamma:C\otimes\mathbf{k}\rightarrow C,\ c\otimes\lambda\mapsto\lambda
c$ and $\gamma^{\prime}:\mathbf{k}\otimes C\rightarrow C,\ \lambda\otimes
c\mapsto\lambda c$ are the canonical $\mathbf{k}$-module isomorphisms. We are
applying both identities to $c$ and using (\ref{pf.lem.cfc.idbar.a.1}) to
expand $\Delta\left(  c\right)  $ on the right hand side.}, we thus have%
\begin{equation}
c=\sum_{i=1}^{m}\epsilon\left(  d_{i}\right)  c_{i}%
\label{pf.lem.cfc.idbar.a.coalg1c}%
\end{equation}
and%
\begin{equation}
c=\sum_{i=1}^{m}\epsilon\left(  c_{i}\right)  d_{i}%
.\label{pf.lem.cfc.idbar.a.coalg2c}%
\end{equation}
Applying the map $\epsilon$ to both sides of (\ref{pf.lem.cfc.idbar.a.coalg2c}%
), we obtain%
\begin{equation}
\epsilon\left(  c\right)  =\epsilon\left(  \sum_{i=1}^{m}\epsilon\left(
c_{i}\right)  d_{i}\right)  =\sum_{i=1}^{m}\epsilon\left(  c_{i}\right)
\epsilon\left(  d_{i}\right)  .\label{pf.lem.cfc.idbar.a.epsc=}%
\end{equation}

Now, applying the map $\idbar\otimes\idbar$ to both sides of the equality
(\ref{pf.lem.cfc.idbar.a.1}), we obtain%
\begin{align*}
&  \left(  \idbar\otimes\idbar\right)  \left(  \Delta\left(  c\right)  \right)
\\
&  =\left(  \idbar\otimes\idbar\right)  \left(  \sum_{i=1}^{m}c_{i}\otimes
d_{i}\right)  =\sum_{i=1}^{m}\underbrace{\idbar\left(  c_{i}\right)
}_{\substack{=c_{i}-\epsilon\left(  c_{i}\right)  1_{C}\\\text{(by the
definition of }\idbar\text{)}}}\otimes\underbrace{\idbar\left(  d_{i}\right)
}_{\substack{=d_{i}-\epsilon\left(  d_{i}\right)  1_{C}\\\text{(by the
definition of }\idbar\text{)}}}\\
&  =\sum_{i=1}^{m}\left(  c_{i}-\epsilon\left(  c_{i}\right)  1_{C}\right)
\otimes\left(  d_{i}-\epsilon\left(  d_{i}\right)  1_{C}\right) \\
&  =\underbrace{\sum_{i=1}^{m}c_{i}\otimes d_{i}}_{\substack{=\Delta\left(
c\right)  \\\text{(by (\ref{pf.lem.cfc.idbar.a.1}))}}}-\underbrace{\sum
_{i=1}^{m}c_{i}\otimes\left(  \epsilon\left(  d_{i}\right)  1_{C}\right)
}_{=\left(  \sum_{i=1}^{m}\epsilon\left(  d_{i}\right)  c_{i}\right)
\otimes1_{C}}-\underbrace{\sum_{i=1}^{m}\left(  \epsilon\left(  c_{i}\right)
1_{C}\right)  \otimes d_{i}}_{=1_{C}\otimes\left(  \sum_{i=1}^{m}%
\epsilon\left(  c_{i}\right)  d_{i}\right)  }+\underbrace{\sum_{i=1}%
^{m}\left(  \epsilon\left(  c_{i}\right)  1_{C}\right)  \otimes\left(
\epsilon\left(  d_{i}\right)  1_{C}\right)  }_{=\left(  \sum_{i=1}^{m}%
\epsilon\left(  c_{i}\right)  \epsilon\left(  d_{i}\right)  \right)
1_{C}\otimes1_{C}}\\
&  =\Delta\left(  c\right)  -\underbrace{\left(  \sum_{i=1}^{m}\epsilon\left(
d_{i}\right)  c_{i}\right)  }_{\substack{=c\\\text{(by
(\ref{pf.lem.cfc.idbar.a.coalg1c}))}}}\otimes1_{C}-1_{C}\otimes
\underbrace{\left(  \sum_{i=1}^{m}\epsilon\left(  c_{i}\right)  d_{i}\right)
}_{\substack{=c\\\text{(by (\ref{pf.lem.cfc.idbar.a.coalg2c}))}}%
}+\underbrace{\left(  \sum_{i=1}^{m}\epsilon\left(  c_{i}\right)
\epsilon\left(  d_{i}\right)  \right)  }_{\substack{=\epsilon\left(  c\right)
\\\text{(by (\ref{pf.lem.cfc.idbar.a.epsc=}))}}}1_{C}\otimes1_{C}\\
&  =\Delta\left(  c\right)  -c\otimes1_{C}-1_{C}\otimes c+\epsilon\left(
c\right)  1_{C}\otimes1_{C}=\delta\left(  c\right)
\end{align*}
(by the definition of $\delta$). Thus,%
\begin{equation}
\delta\left(  c\right)  =\left(  \idbar\otimes\idbar\right)  \left(
\Delta\left(  c\right)  \right)  =\left(  \left(  \idbar\otimes\idbar\right)
\circ\Delta\right)  \left(  c\right)  . \label{pf.lem.cfc.idbar.a.at}%
\end{equation}
Since we have proved this for all $c\in C$, we thus obtain $\delta=\left(
\idbar\otimes\idbar\right)  \circ\Delta$. This proves Lemma
\ref{lem.cfc.idbar} \textbf{(a)}.

\textbf{(b)} Let $n\in\mathbb{N}$. Definition \ref{def.con-fil-coal}
\textbf{(b)} yields $1_{C}\in C_{\leq0}\subseteq C_{\leq n}$ (by
(\ref{eq.def.fil-coal.chain})). Now, for each $c\in C_{\leq n}$, the
definition of $\idbar$ yields%
\[
\idbar\left(  c\right)  =c-\epsilon\left(  c\right)  1_{C}\in C_{\leq
n}\ \ \ \ \ \ \ \ \ \ \left(  \text{since }c\in C_{\leq n}\text{ and }1_{C}\in
C_{\leq n}\right)  .
\]
In other words, we have $\idbar\left(  C_{\leq n}\right)  \subseteq C_{\leq
n}$. This proves Lemma \ref{lem.cfc.idbar} \textbf{(b)}.

\textbf{(c)} The filtered $\mathbf{k}$-coalgebra $C$ is connected. In other
words, the restriction $\epsilon\mid_{C_{\leq0}}$ is a $\mathbf{k}$-module
isomorphism from $C_{\leq0}$ to $\mathbf{k}$ (by Definition
\ref{def.con-fil-coal} \textbf{(a)}). Thus, this restriction $\epsilon
\mid_{C_{\leq0}}$ is injective. Also, Definition \ref{def.con-fil-coal}
\textbf{(b)} yields $1_{C}\in C_{\leq0}$ and $\epsilon\left(  1_{C}\right)
=1_{\mathbf{k}}$.

Now, let $c\in C_{\leq0}$. Set $d=\epsilon\left(  c\right)  1_{C}$. Then,
$d\in C_{\leq0}$ (since $1_{C}\in C_{\leq0}$). From $d=\epsilon\left(
c\right)  1_{C}$, we obtain%
\[
\epsilon\left(  d\right)  =\epsilon\left(  \epsilon\left(  c\right)
1_{C}\right)  =\epsilon\left(  c\right)  \epsilon\left(  1_{C}\right)
=\epsilon\left(  c\right)  \ \ \ \ \ \ \ \ \ \ \left(  \text{since }%
\epsilon\left(  1_{C}\right)  =1_{\mathbf{k}}\right)  .
\]
Since $\epsilon\mid_{C_{\leq0}}$ is injective, this entails $d=c$ (because
both $d$ and $c$ belong to $C_{\leq0}$). Therefore, $c=d=\epsilon\left(
c\right)  1_{C}$. Now, the definition of $\idbar$ yields $\idbar\left(
c\right)  =c-\epsilon\left(  c\right)  1_{C}=0$ (since $c=\epsilon\left(
c\right)  1_{C}$). Since we have proved this for each $c\in C_{\leq0}$, we
thus know that $\idbar\left(  C_{\leq0}\right)  =0$. This proves Lemma
\ref{lem.cfc.idbar} \textbf{(c)}.
\end{proof}

\begin{proof}
[Proof of Proposition \ref{prop.cfc.delta2}.]\textbf{(a)} Define a
$\mathbf{k}$-linear map $\idbar:C\rightarrow C$ as in Lemma
\ref{lem.cfc.idbar}.

Now, let $n>0$ be an integer. Lemma \ref{lem.cfc.idbar} \textbf{(a)} yields
$\delta=\left(  \idbar\otimes\idbar\right)  \circ\Delta$. Thus,%
\begin{align*}
\delta\left(  C_{\leq n}\right)   &  =\left(  \left(  \idbar\otimes
\idbar\right)  \circ\Delta\right)  \left(  C_{\leq n}\right)  =\left(
\idbar\otimes\idbar\right)  \left(  \Delta\left(  C_{\leq n}\right)  \right)
\\
&  \subseteq\left(  \idbar\otimes\idbar\right)  \left(  \sum_{i=0}^{n}C_{\leq
i}\otimes C_{\leq n-i}\right)  \ \ \ \ \ \ \ \ \ \ \left(  \text{by
(\ref{eq.def.fil-coal.Del})}\right)  \\
&  =\sum_{i=0}^{n}\idbar\left(  C_{\leq i}\right)  \otimes\idbar\left(
C_{\leq n-i}\right)  =\sum_{i=1}^{n-1}\underbrace{\idbar\left(  C_{\leq
i}\right)  }_{\substack{\subseteq C_{\leq i}\\\text{(by Lemma
\ref{lem.cfc.idbar} \textbf{(b)})}}}\otimes\underbrace{\idbar\left(  C_{\leq
n-i}\right)  }_{\substack{\subseteq C_{\leq n-i}\\\text{(by Lemma
\ref{lem.cfc.idbar} \textbf{(b)})}}}\\
&  \ \ \ \ \ \ \ \ \ \ \ \ \ \ \ \ \ \ \ \ \left(
\begin{array}
[c]{c}%
\text{here, we have removed the addends for }i=0\text{ and}\\
\text{for }i=n\text{ from the sum, since they are }0\text{ (because}\\
\text{Lemma \ref{lem.cfc.idbar} \textbf{(c)} yields }\idbar\left(  C_{\leq
0}\right)  =0\text{)}%
\end{array}
\right)  \\
&  \subseteq\sum_{i=1}^{n-1}C_{\leq i}\otimes C_{\leq n-i}.
\end{align*}
This proves Proposition \ref{prop.cfc.delta2} \textbf{(a)}.

\textbf{(b)} Let $f:C\rightarrow C$ be a $\mathbf{k}$-coalgebra homomorphism
satisfying $f\left(  1_{C}\right)  =1_{C}$. The definition of a coalgebra
homomorphism thus yields $\left(  f\otimes f\right)  \circ\Delta=\Delta\circ
f$ and $\epsilon=\epsilon\circ f$.

Let $c\in C$. The definition of $\delta$ yields $\delta\left(  c\right)
=\Delta\left(  c\right)  -c\otimes1_{C}-1_{C}\otimes c+\epsilon\left(
c\right)  1_{C}\otimes1_{C}$. Applying the map $f\otimes f$ to both sides of
this equality, we obtain%
\begin{align*}
&  \left(  f\otimes f\right)  \left(  \delta\left(  c\right)  \right) \\
&  =\left(  f\otimes f\right)  \left(  \Delta\left(  c\right)  -c\otimes
1_{C}-1_{C}\otimes c+\epsilon\left(  c\right)  1_{C}\otimes1_{C}\right) \\
&  =\underbrace{\left(  \left(  f\otimes f\right)  \circ\Delta\right)
}_{=\Delta\circ f}\left(  c\right)  -f\left(  c\right)  \otimes
\underbrace{f\left(  1_{C}\right)  }_{=1_{C}}-\underbrace{f\left(
1_{C}\right)  }_{=1_{C}}\otimes f\left(  c\right)  +\underbrace{\epsilon
}_{=\epsilon\circ f}\left(  c\right)  \underbrace{f\left(  1_{C}\right)
}_{=1_{C}}\otimes\underbrace{f\left(  1_{C}\right)  }_{=1_{C}}\\
&  =\Delta\left(  f\left(  c\right)  \right)  -f\left(  c\right)  \otimes
1_{C}-1_{C}\otimes f\left(  c\right)  +\epsilon\left(  f\left(  c\right)
\right)  1_{C}\otimes1_{C}.
\end{align*}
Comparing this with%
\begin{align*}
\left(  \delta\circ f\right)  \left(  c\right)   &  =\delta\left(  f\left(
c\right)  \right)  =\Delta\left(  f\left(  c\right)  \right)  -f\left(
c\right)  \otimes1_{C}-1_{C}\otimes f\left(  c\right)  +\epsilon\left(
f\left(  c\right)  \right)  1_{C}\otimes1_{C}\\
&  \ \ \ \ \ \ \ \ \ \ \left(  \text{by the definition of }\delta\right)  ,
\end{align*}
we obtain $\left(  \delta\circ f\right)  \left(  c\right)  =\left(  f\otimes
f\right)  \left(  \delta\left(  c\right)  \right)  =\left(  \left(  f\otimes
f\right)  \circ\delta\right)  \left(  c\right)  $. Since we have proved this
for each $c\in C$, we thus obtain $\delta\circ f=\left(  f\otimes f\right)
\circ\delta$. This proves Proposition \ref{prop.cfc.delta2} \textbf{(b)}.

\textbf{(c)} Definition \ref{def.con-fil-coal} \textbf{(b)} yields $1_{C}\in
C_{\leq0}$ and $\epsilon\left(  1_{C}\right)  =1_{\mathbf{k}}$.

Let $c\in\left(  \operatorname*{Ker}\delta\right)  \cap\left(
\operatorname*{Ker}\epsilon\right)  $. Thus, $\delta\left(  c\right)  =0$ and
$\epsilon\left(  c\right)  =0$. From $\delta\left(  c\right)  =0$, we obtain%
\begin{align*}
0  &  =\delta\left(  c\right)  =\Delta\left(  c\right)  -c\otimes1_{C}%
-1_{C}\otimes c+\underbrace{\epsilon\left(  c\right)  }_{=0}1_{C}\otimes
1_{C}\ \ \ \ \ \ \ \ \ \ \left(  \text{by the definition of }\delta\right) \\
&  =\Delta\left(  c\right)  -c\otimes1_{C}-1_{C}\otimes c.
\end{align*}
In other words, $\Delta\left(  c\right)  =c\otimes1_{C}+1_{C}\otimes c$. In
other words, the element $c$ of $C$ is primitive. In other words,
$c\in\operatorname*{Prim}C$. Since we have proved this for each $c\in\left(
\operatorname*{Ker}\delta\right)  \cap\left(  \operatorname*{Ker}%
\epsilon\right)  $, we thus obtain $\left(  \operatorname*{Ker}\delta\right)
\cap\left(  \operatorname*{Ker}\epsilon\right)  \subseteq\operatorname*{Prim}%
C$.

Now, let $d\in\operatorname*{Prim}C$. Thus, the element $d$ of $C$ is
primitive, i.e. $\Delta\left(  d\right)  =d\otimes1_{C}%
+1_{C}\otimes d$. Hence, Lemma \ref{lem.coalg.primitive-e0} (applied to
$a=1_{C}$ and $b=1_{C}$) yields $\epsilon\left(  d\right)  =0$ (since
$\epsilon\left(  1_{C}\right)  =1_{\mathbf{k}}=1$). Hence, $d\in
\operatorname*{Ker}\epsilon$.

Furthermore, the definition of $\delta$ yields%
\[
\delta\left(  d\right)  =\underbrace{\Delta\left(  d\right)  -d\otimes
1_{C}-1_{C}\otimes d}_{\substack{=0\\\text{(since }\Delta\left(  d\right)
=d\otimes1_{C}+1_{C}\otimes d\text{)}}}+\underbrace{\epsilon\left(  d\right)
}_{=0}1_{C}\otimes1_{C}=0.
\]
Hence, $d\in\operatorname*{Ker}\delta$. Combining this with $d\in
\operatorname*{Ker}\epsilon$, we obtain $d\in\left(  \operatorname*{Ker}%
\delta\right)  \cap\left(  \operatorname*{Ker}\epsilon\right)  $. Since we
have proved this for each $d\in\operatorname*{Prim}C$, we have thus shown that
$\operatorname*{Prim}C\subseteq\left(  \operatorname*{Ker}\delta\right)
\cap\left(  \operatorname*{Ker}\epsilon\right)  $. Combining this with
$\left(  \operatorname*{Ker}\delta\right)  \cap\left(  \operatorname*{Ker}%
\epsilon\right)  \subseteq\operatorname*{Prim}C$, we obtain the claim of
Proposition \ref{prop.cfc.delta2} \textbf{(c)}.

\textbf{(d)} This follows from Proposition \ref{prop.cfc.delta2} \textbf{(c)},
since both $\operatorname*{Ker}\delta$ and $\operatorname*{Ker}\epsilon$ are
$\mathbf{k}$-submodules of $C$.

\textbf{(e)} Definition \ref{def.con-fil-coal} \textbf{(b)} yields $1_{C}\in
C_{\leq0}\subseteq C_{\leq1}$ (by (\ref{eq.def.fil-coal.chain})). Thus,
\begin{align*}
\delta\left(  1_{C}\right)   &  \in\delta\left(  C_{\leq1}\right)
\subseteq\sum_{i=1}^{1-1}C_{\leq i}\otimes C_{\leq1-i}%
\ \ \ \ \ \ \ \ \ \ \left(  \text{by Proposition \ref{prop.cfc.delta2}
\textbf{(a)}, applied to }n=1\right) \\
&  =\left(  \text{empty sum}\right)  =0.
\end{align*}
In other words, $\delta\left(  1_{C}\right)  =0$. Hence, $1_{C}\in
\operatorname*{Ker}\delta$. Thus, $\mathbf{k}\cdot1_{C}\subseteq
\operatorname*{Ker}\delta$.

Also, Proposition \ref{prop.cfc.delta2} \textbf{(c)} yields
\[
\operatorname*{Prim}C=\left(  \operatorname*{Ker}\delta\right)  \cap\left(
\operatorname*{Ker}\epsilon\right)  \subseteq\operatorname*{Ker}\delta.
\]
Hence,%
\begin{equation}
\underbrace{\mathbf{k}\cdot1_{C}}_{\subseteq\operatorname*{Ker}\delta
}+\underbrace{\operatorname*{Prim}C}_{\subseteq\operatorname*{Ker}\delta
}\subseteq\operatorname*{Ker}\delta+\operatorname*{Ker}\delta\subseteq
\operatorname*{Ker}\delta. \label{pf.prop.cfc.delta2.e.3}%
\end{equation}

Definition \ref{def.con-fil-coal} \textbf{(b)} yields $1_{C}\in C_{\leq0}$ and
$\epsilon\left(  1_{C}\right)  =1_{\mathbf{k}}$.

Let $u\in\operatorname*{Ker}\delta$. Thus, $\delta\left(  u\right)  =0$. Set
$v=u-\epsilon\left(  u\right)  1_{C}$. Thus,%
\[
\delta\left(  v\right)  =\delta\left(  u-\epsilon\left(  u\right)
1_{C}\right)  =\underbrace{\delta\left(  u\right)  }_{=0}-\epsilon\left(
u\right)  \underbrace{\delta\left(  1_{C}\right)  }_{=0}=0-\epsilon\left(
u\right)  0=0,
\]
so that $v\in\operatorname*{Ker}\delta$. Furthermore, from $v=u-\epsilon
\left(  u\right)  1_{C}$, we obtain
\[
\epsilon\left(  v\right)  =\epsilon\left(  u-\epsilon\left(  u\right)
1_{C}\right)  =\epsilon\left(  u\right)  -\epsilon\left(  u\right)
\underbrace{\epsilon\left(  1_{C}\right)  }_{=1_{\mathbf{k}}}=\epsilon\left(
u\right)  -\epsilon\left(  u\right)  =0,
\]
so that $v\in\operatorname*{Ker}\epsilon$. Combining this with $v\in
\operatorname*{Ker}\delta$, we obtain $v\in\left(  \operatorname*{Ker}%
\delta\right)  \cap\left(  \operatorname*{Ker}\epsilon\right)
=\operatorname*{Prim}C$ (by Proposition \ref{prop.cfc.delta2} \textbf{(c)}).
Now, from $v=u-\epsilon\left(  u\right)  1_{C}$, we obtain%
\[
u=\underbrace{\epsilon\left(  u\right)  }_{\in\mathbf{k}}1_{C}+\underbrace{v}%
_{\in\operatorname*{Prim}C}\in\mathbf{k}\cdot1_{C}+\operatorname*{Prim}C.
\]

We thus have shown that $u\in\mathbf{k}\cdot1_{C}+\operatorname*{Prim}C$ for
each $u\in\operatorname*{Ker}\delta$. In other words,%
\[
\operatorname*{Ker}\delta\subseteq\mathbf{k}\cdot1_{C}+\operatorname*{Prim}C.
\]
Combining this with (\ref{pf.prop.cfc.delta2.e.3}), we obtain
$\operatorname*{Ker}\delta=\mathbf{k}\cdot1_{C}+\operatorname*{Prim}C$ and prove Proposition~\ref{prop.cfc.delta2}~\textbf{(e)}.
\end{proof}

\subsection{\label{sec.pf.cor.id-f.cfc}Proofs of the corollaries from Section
\ref{sec.cfc}}

\begin{proof}
[Proof of Corollary \ref{cor.id-f.cfc}.]From $e\left(  1_{C}\right)  =1_{C}$
and $f\left(  1_{C}\right)  =1_{C}$, we obtain
\[
    \left(  e-f\right)  \left( 1_{C}\right)  =1_{C}-1_{C}=0.
\]
Hence, $1_{C}\in\operatorname*{Ker}\left(
e-f\right)  $, so that $\mathbf{k}\cdot1_{C}\subseteq\operatorname*{Ker}%
\left(  e-f\right)  $.

Define the $\mathbf{k}$-linear map $\delta:C\rightarrow C\otimes C$ as in
Proposition \ref{prop.cfc.delta2}. Proposition \ref{prop.cfc.delta2}
\textbf{(b)} yields that $\left(  f\otimes f\right)  \circ\delta=\delta\circ
f$. Likewise, $\left(  e\otimes e\right)  \circ\delta=\delta\circ e$.
Moreover, Proposition \ref{prop.cfc.delta2} \textbf{(e)} yields
\[
\operatorname*{Ker}\delta=\underbrace{\mathbf{k}\cdot1_{C}}_{\subseteq
\operatorname*{Ker}\left(  e-f\right)  }+\underbrace{\operatorname*{Prim}%
C}_{\substack{\subseteq\operatorname*{Ker}\left(  e-f\right)  \\\text{(by
(\ref{eq.cor.id-f.cfc.ass-Ker}))}}}\subseteq\operatorname*{Ker}\left(
e-f\right)  +\operatorname*{Ker}\left(  e-f\right)  \subseteq
\operatorname*{Ker}\left(  e-f\right)  .
\]

However, Proposition \ref{prop.cfc.delta2} \textbf{(a)} shows that%
\[
\delta\left(  C_{\leq n}\right)  \subseteq\sum_{i=1}^{n-1}C_{\leq i}\otimes
C_{\leq n-i}\ \ \ \ \ \ \ \ \ \ \text{for each }n>p
\]
(and, more generally, for each $n>0$, but we only need this in the case $n>p$).

Moreover, (\ref{eq.def.fil-coal.chain}) yields $C_{\leq1}+C_{\leq2}%
+\cdots+C_{\leq p}\subseteq C_{\leq p}$. Therefore,%
\[
\left(  e-f\right)  \left(  C_{\leq1}+C_{\leq2}+\cdots+C_{\leq p}\right)
\subseteq\left(  e-f\right)  \left(  C_{\leq p}\right)  =0
\]
(by (\ref{eq.cor.id-f.cfc.ass-ann})), so that $\left(  e-f\right)  \left(
C_{\leq1}+C_{\leq2}+\cdots+C_{\leq p}\right)  =0$.

Hence, Theorem \ref{thm.id-f.gen} (applied to $D=C$ and $D_{i}=C_{\leq i}$)
shows that for any integer $u>p$, we have%
\begin{equation}
\left(  e-f\right)  ^{u-p}\left(  C_{\leq u}\right)  \subseteq
\operatorname*{Ker}\delta\label{pf.cor.id-f.cfc.claim-1a}%
\end{equation}
and%
\begin{equation}
\left(  e-f\right)  ^{u-p+1}\left(  C_{\leq u}\right)  =0.
\label{pf.cor.id-f.cfc.claim-1b}%
\end{equation}

We are now close to proving Corollary \ref{cor.id-f.cfc}. Let us begin with
part \textbf{(a)}:

\textbf{(a)} The map $f$ is a $\mathbf{k}$-coalgebra homomorphism, and thus
satisfies $\epsilon\circ f=\epsilon$ (by the definition of a $\mathbf{k}%
$-coalgebra homomorphism). Similarly, $\epsilon\circ e=\epsilon$. Since the
map $\epsilon$ is $\mathbf{k}$-linear, we have
\[
\epsilon\circ\left(  e-f\right)  =\underbrace{\epsilon\circ e}_{=\epsilon
}-\underbrace{\epsilon\circ f}_{=\epsilon}=\epsilon-\epsilon=0.
\]

Now, let $u>p$ be an integer. Thus, $\left(  e-f\right)  ^{u-p}=\left(
e-f\right)  \circ\left(  e-f\right)  ^{u-p-1}$. Hence,%
\[
\epsilon\circ\left(  e-f\right)  ^{u-p}=\underbrace{\epsilon\circ\left(
e-f\right)  }_{=0}\circ\left(  e-f\right)  ^{u-p-1}=0\circ\left(  e-f\right)
^{u-p-1}=0.
\]
Therefore, $\left(  e-f\right)  ^{u-p}\left(  C_{\leq u}\right)
\subseteq\operatorname*{Ker}\epsilon$. Combining this with
(\ref{pf.cor.id-f.cfc.claim-1a}), we find, by
Proposition~\ref{prop.cfc.delta2}~\textbf{(c)}, that 
\[
\left(  e-f\right)  ^{u-p}\left(
C_{\leq u}\right)  \subseteq\left(  \operatorname*{Ker}\delta\right)
\cap\left(  \operatorname*{Ker}\epsilon\right)  =\operatorname*{Prim}C.
\]
This proves Corollary
\ref{cor.id-f.cfc} \textbf{(a)}.

\textbf{(b)} Let $u\geq p$ be an integer. We must prove that $\left(
e-f\right)  ^{u-p+1}\left(  C_{\leq u}\right)  =0$. If $u>p$, then this
follows from (\ref{pf.cor.id-f.cfc.claim-1b}). Thus, for the rest of this
proof, we assume (without loss of generality) that we don't have $u>p$. Hence, $u=p$ (since we have
$u\geq p$). Thus, $\left(  e-f\right)  ^{u-p+1}=\left(  e-f\right)
^{p-p+1}=\left(  e-f\right)  ^{1}=e-f$ and $C_{\leq u}=C_{\leq p}$, so that%
\[
\left(  e-f\right)  ^{u-p+1}\left(  C_{\leq u}\right)  =\left(  e-f\right)
\left(  C_{\leq p}\right)  =0
\]
(by (\ref{eq.cor.id-f.cfc.ass-ann})). This proves Corollary \ref{cor.id-f.cfc}
\textbf{(b)}.
\end{proof}

\begin{proof}
[Proof of Corollary \ref{cor.id-f.cfc1}.]Define the $\mathbf{k}$-linear map
$\delta:C\rightarrow C\otimes C$ as in Proposition \ref{prop.cfc.delta2}. Just
as we did in the proof of Corollary \ref{cor.id-f.cfc}, we can show that
$\operatorname*{Ker}\delta\subseteq\operatorname*{Ker}\left(  e-f\right)  $.
As we have seen in the proof of Proposition \ref{prop.cfc.delta2}
\textbf{(e)}, we have $\delta\left(  C_{\leq1}\right)  \subseteq0$. Hence,
\[
    C_{\leq1}
    \subseteq\operatorname*{Ker}\delta
    \subseteq\operatorname*{Ker}\left(  e-f\right).
\]
In other words,
\begin{equation}
\left(  e-f\right)  \left(  C_{\leq1}\right)  =0.
\label{pf.cor.id-f.cfc1.ec=0}%
\end{equation}
Hence, we can apply Corollary \ref{cor.id-f.cfc} to $p=1$. In particular,
applying Corollary \ref{cor.id-f.cfc} \textbf{(a)} to $p=1$, we immediately
obtain Corollary \ref{cor.id-f.cfc1} \textbf{(a)}, whereas applying Corollary
\ref{cor.id-f.cfc} \textbf{(b)} to $p=1$ yields Corollary \ref{cor.id-f.cfc1}
\textbf{(b)}.
\end{proof}

\begin{proof}
[Proof of Corollary \ref{cor.id-f.cfc1id}.]Clearly, $\operatorname*{id}%
:C\rightarrow C$ is a $\mathbf{k}$-coalgebra homomorphism such that
$\operatorname*{id}\left(  1_{C}\right)  =1_{C}$. Furthermore, $f\circ
\operatorname*{id}=f=\operatorname*{id}\circ f$. Hence, we can apply Corollary
\ref{cor.id-f.cfc1} to $e=\operatorname*{id}$. As a result, we obtain
precisely the claims of Corollary \ref{cor.id-f.cfc1id}.
\end{proof}

\subsection{\label{sec.pf.filhopf}Proofs for Section \ref{sec.filhopf}}

Before we prove the claims left unproved in Section \ref{sec.filhopf}, let us
recall the defining property of the antipode of a Hopf algebra (see, e.g.,
\cite[(1.4.3)]{GriRei}):

\begin{remark}
\label{rmk.hopf.antipode-pro}Let $H$ be a $\mathbf{k}$-Hopf algebra with
antipode $S$. Let $1_{H}$ denote the unity of the $\mathbf{k}$-algebra $H$.
Let $m:H\otimes H\rightarrow H$ be the $\mathbf{k}$-linear map that sends each
pure tensor $x\otimes y\in H\otimes H$ to the product $xy\in H$. Let
$u:\mathbf{k}\rightarrow H$ be the $\mathbf{k}$-linear map that sends
$1_{\mathbf{k}}$ to $1_{H}$. Then,%
\begin{align}
m\circ\left(  S\otimes\operatorname*{id}\nolimits_{H}\right)  \circ\Delta &
=u\circ\epsilon\ \ \ \ \ \ \ \ \ \ \text{and}%
\label{eq.rmk.hopf.antipode-pro.1}\\
m\circ\left(  \operatorname*{id}\nolimits_{H}\otimes S\right)  \circ\Delta &
=u\circ\epsilon. \label{eq.rmk.hopf.antipode-pro.2}%
\end{align}

\end{remark}

\begin{proof}
[Proof of Corollary \ref{cor.id-S2.p}.]Lemma \ref{lem.bialg.antip-props}
\textbf{(b)} yields $S\left(  1_{H}\right)  =1_{H}$. Hence, it easily follows
that $S^{2}\left(  1_{H}\right)  =1_{H}$. Moreover, Lemma
\ref{lem.bialg.antip-props} \textbf{(a)} yields that the map $S^{2}%
:H\rightarrow H$ is a $\mathbf{k}$-coalgebra homomorphism. Of course, the map
$\operatorname*{id}:H\rightarrow H$ is a $\mathbf{k}$-coalgebra homomorphism
as well, and satisfies $\operatorname*{id}\left(  1_{H}\right)  =1_{H}$.
Furthermore, every $x\in\operatorname*{Prim}H$ is a primitive element of $H$
and therefore satisfies $x\in\operatorname*{Ker}\left(  \operatorname*{id}%
-S^{2}\right)  $ (since Lemma \ref{lem.bialg.antip-props} \textbf{(d)} yields
$S^{2}\left(  x\right)  =x$). Hence, we have $\operatorname*{Prim}%
H\subseteq\operatorname*{Ker}\left(  \operatorname*{id}-S^{2}\right)  $.
Moreover, $S^{2}\circ\operatorname*{id}=S^{2}=\operatorname*{id}\circ S^{2}$.
Furthermore, $p$ is a positive integer and satisfies $\left(
\operatorname*{id}-S^{2}\right)  \left(  H_{\leq p}\right)  =0$ (by
(\ref{eq.cor.id-S2.p.ass-ann})). Hence, we can apply Corollary
\ref{cor.id-f.cfc} to $C=H$ and $C_{\leq i}=H_{\leq i}$ and
$e=\operatorname*{id}$ and $f=S^{2}$. Doing so, we immediately obtain

\begin{itemize}
\item that (\ref{eq.cor.id-S2.p.claim1}) holds for any integer $u>p$ (by
applying Corollary \ref{cor.id-f.cfc} \textbf{(a)}), and

\item that (\ref{eq.cor.id-S2.p.claim2}) holds for any integer $u\geq p$ (by
applying Corollary \ref{cor.id-f.cfc} \textbf{(b)}).
\end{itemize}

\noindent It thus remains to prove that (\ref{eq.cor.id-S2.p.claim3}) holds
for any integer $u>p$. So let us do this now.

First, we shall show that $\left(  \operatorname*{id}+S\right)  \left(
\operatorname*{Prim}H\right)  =0$.

Indeed, each $x\in\operatorname*{Prim}H$ is a primitive element of $H$ and
therefore satisfies%
\[
\left(  \operatorname*{id}+S\right)  \left(  x\right)  =x+S\left(  x\right)
=0\ \ \ \ \ \ \ \ \ \ \left(  \text{since Lemma \ref{lem.bialg.antip-props}
\textbf{(c)} yields }S\left(  x\right)  =-x\right)  .
\]
In other words, we have $\left(  \operatorname*{id}+S\right)  \left(
\operatorname*{Prim}H\right)  =0$.

Now, let $u>p$ be an integer. Applying the map $\operatorname*{id}+S$ to
(\ref{eq.cor.id-S2.p.claim1}), we obtain%
\[
\left(  \left(  \operatorname*{id}+S\right)  \circ\left(  \operatorname*{id}%
-S^{2}\right)  ^{u-p}\right)  \left(  H_{\leq u}\right)  \subseteq\left(
\operatorname*{id}+S\right)  \left(  \operatorname*{Prim}H\right)  =0.
\]
Thus, (\ref{eq.cor.id-S2.p.claim3}) follows. This completes the proof of
Corollary \ref{cor.id-S2.p}.
\end{proof}

\begin{proof}
[Proof of Corollary \ref{cor.id-S2.1}.]In our above proof of Corollary
\ref{cor.id-S2.p}, we have already shown that

\begin{itemize}
\item we have $S^{2}\left(  1_{H}\right)  =1_{H}$;

\item the map $S^{2}:H\rightarrow H$ is a $\mathbf{k}$-coalgebra homomorphism;

\item we have $\operatorname*{Prim}H\subseteq\operatorname*{Ker}\left(
\operatorname*{id}-S^{2}\right)  $.
\end{itemize}

Hence, we can apply Corollary \ref{cor.id-f.cfc1id} to $C=H$ and $C_{\leq
i}=H_{\leq i}$ and $f=S^{2}$. Doing so, we immediately obtain

\begin{itemize}
\item that (\ref{eq.cor.id-S2.1.claim1}) holds for any integer $u>1$ (by
applying Corollary \ref{cor.id-f.cfc1id} \textbf{(a)}), and

\item that (\ref{eq.cor.id-S2.1.claim2}) holds for any positive integer $u$
(by applying Corollary \ref{cor.id-f.cfc1id} \textbf{(b)}).
\end{itemize}

\noindent It thus remains to prove that (\ref{eq.cor.id-S2.1.claim3}) holds
for any integer $u>1$. But this can be deduced from
(\ref{eq.cor.id-S2.1.claim1}) in the same way as we deduced
(\ref{eq.cor.id-S2.p.claim3}) from (\ref{eq.cor.id-S2.p.claim1}) in our above
proof of Corollary \ref{cor.id-S2.p}. Thus, the proof of Corollary
\ref{cor.id-S2.1} is complete.
\end{proof}

\subsection{\label{sec.pf.hopf}Proofs for Section \ref{sec.hopf}}

We shall next focus on proving the claims left unproven in Section
\ref{sec.hopf}. Before we do so, let us first collect a few basic properties
of connected graded Hopf algebras into a lemma for convenience:

\begin{lemma}
\label{lem.cghopf.basics}Let $H$ be a connected graded $\mathbf{k}$-Hopf
algebra with unity $1_{H}$ and antipode $S$. Then:

\textbf{(a)} If $n$ is a positive integer, and if $x$ is an element of $H_{n}%
$, then we have%
\[
\Delta\left(  x\right)  =1_{H}\otimes x+x\otimes1_{H}%
+w\ \ \ \ \ \ \ \ \ \ \text{for some }w\in\sum_{k=1}^{n-1}H_{k}\otimes
H_{n-k}.
\]

\textbf{(b)} We have $H_{1}\subseteq\operatorname*{Prim}H$.

\textbf{(c)} We have $S\left(  ab\right)  =ba$ for any $a,b\in H_{1}$.
\end{lemma}

\begin{proof}
[Proof of Lemma \ref{lem.cghopf.basics}.]\textbf{(a)} This is well-known; see
\cite[Exercise 1.3.20 (h)]{GriRei} or \cite[Proposition II.1.1]{Mancho06} or
\cite[Theorem 2.18]{Preiss16} for a proof.

\textbf{(b)} We need to show that each $x\in H_{1}$ is primitive, i.e.,
satisfies $\Delta\left(  x\right)  =x\otimes1_{H}+1_{H}\otimes x$. But this
follows easily by applying Lemma \ref{lem.cghopf.basics} \textbf{(a)} to $n=1$
(and observing that the sum $\sum_{k=1}^{n-1}H_{k}\otimes H_{n-k}$ is empty
for $n=1$).

\textbf{(c)} Let $a,b\in H_{1}$. Then, $a\in H_{1}\subseteq
\operatorname*{Prim}H$ (by Lemma \ref{lem.cghopf.basics} \textbf{(b)}). In
other words, the element $a$ of $H$ is primitive. Hence, $S\left(  a\right)
=-a$ (by Lemma \ref{lem.bialg.antip-props} \textbf{(c)}, applied to $x=a$).
Similarly, $S\left(  b\right)  =-b$. However, it is well-known (see, e.g.,
\cite[Proposition 1.2.17 (1)]{HecSch20} or \cite[Proposition 1.4.10]{GriRei}
or \cite[Proposition 7.1.9 (a)]{Radfor12}) that the antipode $S$ of $H$ is a
$\mathbf{k}$-algebra anti-endomorphism, i.e., that it satisfies $S\left(
1_{H}\right)  =1_{H}$ and
\begin{equation}
S\left(  uv\right)  =S\left(  v\right)  S\left(  u\right)
\ \ \ \ \ \ \ \ \ \ \text{for all }u,v\in H. \label{pf.lem.cghopf.basics.c.uv}%
\end{equation}
Applying (\ref{pf.lem.cghopf.basics.c.uv}) to $u=a$ and $v=b$, we obtain
$S\left(  ab\right)  =\underbrace{S\left(  b\right)  }_{=-b}%
\underbrace{S\left(  a\right)  }_{=-a}=\left(  -b\right)  \left(  -a\right)
=ba$. This proves Lemma \ref{lem.cghopf.basics} \textbf{(c)}.
\end{proof}

\begin{proof}
[Proof of Corollary \ref{cor.id-S2.gr1}.]As we know, the graded $\mathbf{k}%
$-Hopf algebra $H$ automatically becomes a filtered $\mathbf{k}$-Hopf algebra
with filtration $\left(  H_{\leq0},H_{\leq1},H_{\leq2},\ldots\right)  $
defined by setting%
\[
H_{\leq n}:=\bigoplus_{i=0}^{n}H_{i}\ \ \ \ \ \ \ \ \ \ \text{for all }%
n\in\mathbb{N}.
\]
This filtered $\mathbf{k}$-Hopf algebra $H$ is connected, since $H_{\leq
0}=H_{0}$. Thus, Corollary \ref{cor.id-S2.1} can be applied.

Let $u$ be a positive integer. Then, the definition of $H_{\leq u}$ yields
$H_{u}\subseteq H_{\leq u}$.

Now, we must prove the three relations (\ref{eq.cor.id-S2.gr1.claim1}),
(\ref{eq.cor.id-S2.gr1.claim3}) and (\ref{eq.cor.id-S2.gr1.claim2}). The third
one is the easiest: From $H_{u}\subseteq H_{\leq u}$, we obtain%
\[
\left(  \operatorname*{id}-S^{2}\right)  ^{u}\left(  H_{u}\right)
\subseteq\left(  \operatorname*{id}-S^{2}\right)  ^{u}\left(  H_{\leq
u}\right)  =0
\]
(by Corollary \ref{cor.id-S2.1} \textbf{(b)}) and therefore $\left(
\operatorname*{id}-S^{2}\right)  ^{u}\left(  H_{u}\right)  =0$. This proves
(\ref{eq.cor.id-S2.gr1.claim2}).

We shall now focus on proving (\ref{eq.cor.id-S2.gr1.claim1}). Indeed, if
$u>1$, then (\ref{eq.cor.id-S2.gr1.claim1}) follows from%
\begin{align*}
\left(  \operatorname*{id}-S^{2}\right)  ^{u-1}\left(  H_{u}\right)   &
\subseteq\left(  \operatorname*{id}-S^{2}\right)  ^{u-1}\left(  H_{\leq
u}\right)  \ \ \ \ \ \ \ \ \ \ \left(  \text{since }H_{u}\subseteq H_{\leq
u}\right) \\
&  \subseteq\operatorname*{Prim}H\ \ \ \ \ \ \ \ \ \ \left(  \text{by
(\ref{eq.cor.id-S2.1.claim1}), since }u>1\right)  .
\end{align*}
Thus, in order to complete the proof of (\ref{eq.cor.id-S2.gr1.claim1}), we
only need to prove it for $u=1$. In other words, we need to prove that
$\left(  \operatorname*{id}-S^{2}\right)  ^{0}\left(  H_{1}\right)
\subseteq\operatorname*{Prim}H$. But this follows from Lemma
\ref{lem.cghopf.basics} \textbf{(b)}, since $\left(  \operatorname*{id}%
-S^{2}\right)  ^{0}=\operatorname*{id}$. This completes our proof of
(\ref{eq.cor.id-S2.gr1.claim1}).

Now, it remains to prove (\ref{eq.cor.id-S2.gr1.claim3}). But we can deduce
(\ref{eq.cor.id-S2.gr1.claim3}) from (\ref{eq.cor.id-S2.gr1.claim1}) in the
same way as we deduced (\ref{eq.cor.id-S2.p.claim3}) from
(\ref{eq.cor.id-S2.p.claim1}) in our above proof of Corollary
\ref{cor.id-S2.p}. This completes the proof of Corollary \ref{cor.id-S2.gr1}.
\end{proof}

\begin{proof}
[Proof of Corollary \ref{cor.id-S2.grp}.]Let $1_{H}$ denote the unity of the
$\mathbf{k}$-algebra $H$.

As in the proof of Corollary \ref{cor.id-S2.gr1}, we know that the graded
$\mathbf{k}$-Hopf algebra $H$ automatically becomes a filtered $\mathbf{k}%
$-Hopf algebra, and this filtered $\mathbf{k}$-Hopf algebra $H$ is connected.

Now, we shall show that
\begin{equation}
\left(  \operatorname*{id}-S^{2}\right)  \left(  H_{\leq p}\right)  =0.
\label{pf.cor.id-S2.grp.1}%
\end{equation}

[\textit{Proof of (\ref{pf.cor.id-S2.grp.1}):} In view of $H_{\leq
p}=\bigoplus_{i=0}^{p}H_{i}$, it will suffice to show that $\left(
\operatorname*{id}-S^{2}\right)  \left(  H_{i}\right)  =0$ for each
$i\in\left\{  0,1,\ldots,p\right\}  $. But (\ref{eq.cor.id-S2.grp.ass-ann})
shows that this is true for all $i\in\left\{  2,3,\ldots,p\right\}  $. Thus,
it remains to show that $\left(  \operatorname*{id}-S^{2}\right)  \left(
H_{i}\right)  =0$ for $i=0$ and for $i=1$.

This is not hard: Lemma \ref{lem.bialg.antip-props} \textbf{(b)} easily yields
$\left(  \operatorname*{id}-S^{2}\right)  \left(  1_{H}\right)  =0$. However,
since the graded Hopf algebra $H$ is connected, it is easy to see that
$H_{0}=\mathbf{k}\cdot1_{H}$. Hence, from $\left(  \operatorname*{id}%
-S^{2}\right)  \left(  1_{H}\right)  =0$, we obtain $\left(
\operatorname*{id}-S^{2}\right)  \left(  H_{0}\right)  =0$. Thus, $\left(
\operatorname*{id}-S^{2}\right)  \left(  H_{i}\right)  =0$ holds for $i=0$.

Each element $x\in H_{1}$ is primitive (by Lemma \ref{lem.cghopf.basics}
\textbf{(b)}) and thus satisfies $S^{2}\left(  x\right)  =x$ (by Lemma
\ref{lem.bialg.antip-props} \textbf{(d)}), so that $\left(  \operatorname*{id}%
-S^{2}\right)  \left(  x\right)  =0$. In other words, $\left(
\operatorname*{id}-S^{2}\right)  \left(  H_{1}\right)  =0$. Thus, $\left(
\operatorname*{id}-S^{2}\right)  \left(  H_{i}\right)  =0$ holds for $i=1$.
Our proof of (\ref{pf.cor.id-S2.grp.1}) is now complete.]

Hence, we can apply Corollary \ref{cor.id-S2.p}. As a result, we obtain
precisely the claims of Corollary \ref{cor.id-S2.grp}.
\end{proof}

\begin{proof}
[Proof of Corollary \ref{cor.id-S2.gr2}.]\textbf{(a)} Let $1_{H}$ denote the
unity of the $\mathbf{k}$-algebra $H$. Define the maps $m$ and $u$ as in
Remark \ref{rmk.hopf.antipode-pro}.

Let $x\in H_{2}$. Then, Lemma \ref{lem.cghopf.basics} \textbf{(a)} (applied to
$n=2$) yields that we have%
\[
\Delta\left(  x\right)  =1_{H}\otimes x+x\otimes1_{H}%
+w\ \ \ \ \ \ \ \ \ \ \text{for some }w\in H_{1}\otimes H_{1}%
\]
(since $\sum_{k=1}^{2-1}H_{k}\otimes H_{2-k}=H_{1}\otimes H_{1}$). Consider
this $w$.

Now, $w$ is a tensor in $H_{1}\otimes H_{1}$. Write this tensor in the form%
\begin{equation}
w=\sum_{i=1}^{k}a_{i}\otimes b_{i} \label{pf.cor.id-S2.gr2.a.3}%
\end{equation}
for some $k\in\mathbb{N}$, some $a_{1},a_{2},\ldots,a_{k}\in H_{1}$ and some
$b_{1},b_{2},\ldots,b_{k}\in H_{1}$.

The elements $a_{1},a_{2},\ldots,a_{k}$ and $b_{1},b_{2},\ldots,b_{k}$ belong
to $H_{1}$ and therefore are primitive (since Lemma \ref{lem.cghopf.basics}
\textbf{(b)} yields $H_{1}\subseteq\operatorname*{Prim}H$). Hence, each
$i\in\left\{  1,2,\ldots,k\right\}  $ satisfies%
\begin{equation}
S\left(  a_{i}\right)  =-a_{i} \label{pf.cor.id-S2.gr2.a.4}%
\end{equation}
(by Lemma \ref{lem.bialg.antip-props} \textbf{(c)}, applied to $a_{i}$ instead
of $x$) and%
\begin{align}
S\left(  a_{i}b_{i}\right)   &  =b_{i}a_{i}\ \ \ \ \ \ \ \ \ \ \left(
\text{by Lemma \ref{lem.cghopf.basics} \textbf{(c)}, applied to }a=a_{i}\text{
and }b=b_{i}\right) \nonumber\\
&  =a_{i}b_{i} \label{pf.cor.id-S2.gr2.a.Saibi}%
\end{align}
(by (\ref{eq.cor.id-S2.gr2.commH1}), applied to $a=b_{i}$ and $b=a_{i}$).

Applying the map $S\otimes\operatorname*{id}\nolimits_{H}:H\otimes
H\rightarrow H\otimes H$ to the equality (\ref{pf.cor.id-S2.gr2.a.3}), we
obtain%
\begin{align}
\left(  S\otimes\operatorname*{id}\nolimits_{H}\right)  \left(  w\right)   &
=\left(  S\otimes\operatorname*{id}\nolimits_{H}\right)  \left(  \sum
_{i=1}^{k}a_{i}\otimes b_{i}\right)  =\sum_{i=1}^{k}\underbrace{S\left(
a_{i}\right)  }_{\substack{=-a_{i}\\\text{(by (\ref{pf.cor.id-S2.gr2.a.4}))}%
}}\otimes\underbrace{\operatorname*{id}\nolimits_{H}\left(  b_{i}\right)
}_{=b_{i}}=-\underbrace{\sum_{i=1}^{k}a_{i}\otimes b_{i}}%
_{\substack{=w\\\text{(by (\ref{pf.cor.id-S2.gr2.a.3}))}}}\nonumber\\
&  =-w. \label{pf.cor.id-S2.gr2.a.5}%
\end{align}

The Hopf algebra $H$ is graded. Hence, its counit $\epsilon$ is a graded map
from $H$ to $\mathbf{k}$. In other words, $\epsilon\left(  H_{i}\right)
\subseteq\mathbf{k}_{i}$ for each $i\in\mathbb{N}$. Thus, $\epsilon\left(
H_{2}\right)  \subseteq\mathbf{k}_{2}=0$ (since the graded $\mathbf{k}$-module
$\mathbf{k}$ is concentrated in degree $0$). Therefore, $\epsilon\left(
x\right)  =0$ (since $x\in H_{2}$).

Lemma \ref{lem.bialg.antip-props} \textbf{(b)} yields $S\left(  1_{H}\right)
=1_{H}$.

Applying both sides of the equality (\ref{eq.rmk.hopf.antipode-pro.1}) to $x$,
we obtain%
\[
\left(  m\circ\left(  S\otimes\operatorname*{id}\nolimits_{H}\right)
\circ\Delta\right)  \left(  x\right)  =\left(  u\circ\epsilon\right)  \left(
x\right)  =u\left(  \underbrace{\epsilon\left(  x\right)  }_{=0}\right)
=u\left(  0\right)  =0.
\]
Therefore,%
\begin{align*}
0  &  =\left(  m\circ\left(  S\otimes\operatorname*{id}\nolimits_{H}\right)
\circ\Delta\right)  \left(  x\right) \\
&  =m\left(  \underbrace{S\left(  1_{H}\right)  }_{=1_{H}}\otimes
\underbrace{\operatorname*{id}\nolimits_{H}\left(  x\right)  }_{=x}+S\left(
x\right)  \otimes\underbrace{\operatorname*{id}\nolimits_{H}\left(
1_{H}\right)  }_{=1_{H}}+\underbrace{\left(  S\otimes\operatorname*{id}%
\nolimits_{H}\right)  \left(  w\right)  }_{\substack{=-w\\\text{(by
(\ref{pf.cor.id-S2.gr2.a.5}))}}}\right) \\
&  \ \ \ \ \ \ \ \ \ \ \ \ \ \ \ \ \ \ \ \ \left(  \text{since }\Delta\left(
x\right)  =1_{H}\otimes x+x\otimes1_{H}+w\right) \\
&  =m\left(  1_{H}\otimes x+S\left(  x\right)  \otimes1_{H}-w\right) \\
&  =1_{H}x+S\left(  x\right)  \cdot1_{H}-m\left(  w\right)
\ \ \ \ \ \ \ \ \ \ \left(  \text{by the definition of }m\right) \\
&  =x+S\left(  x\right)  -m\left(  w\right)  .
\end{align*}
Solving this equality for $S\left(  x\right)  $, we obtain%
\begin{equation}
S\left(  x\right)  =m\left(  w\right)  -x. \label{pf.cor.id-S2.gr2.a.Sx=}%
\end{equation}

Applying the map $m:H\otimes H\rightarrow H$ to the equality
(\ref{pf.cor.id-S2.gr2.a.3}), we obtain%
\begin{equation}
m\left(  w\right)  =m\left(  \sum_{i=1}^{k}a_{i}\otimes b_{i}\right)
=\sum_{i=1}^{k}a_{i}b_{i} \label{pf.cor.id-S2.gr2.a.mw=sum}%
\end{equation}
(by the definition of $m$). Applying the map $S$ to this equality, we obtain%
\begin{align}
S\left(  m\left(  w\right)  \right)   &  =S\left(  \sum_{i=1}^{k}a_{i}%
b_{i}\right)  =\sum_{i=1}^{k}\underbrace{S\left(  a_{i}b_{i}\right)
}_{\substack{=a_{i}b_{i}\\\text{(by (\ref{pf.cor.id-S2.gr2.a.Saibi}))}}%
}=\sum_{i=1}^{k}a_{i}b_{i}\nonumber\\
&  =m\left(  w\right)  \ \ \ \ \ \ \ \ \ \ \left(  \text{by
(\ref{pf.cor.id-S2.gr2.a.mw=sum})}\right)  . \label{pf.cor.id-S2.gr2.a.Smw=}%
\end{align}

Now, applying the map $S$ to both sides of the equality
(\ref{pf.cor.id-S2.gr2.a.Sx=}), we obtain%
\begin{align*}
S\left(  S\left(  x\right)  \right)   &  =S\left(  m\left(  w\right)
-x\right)  =\underbrace{S\left(  m\left(  w\right)  \right)  }%
_{\substack{=m\left(  w\right)  \\\text{(by (\ref{pf.cor.id-S2.gr2.a.Smw=}))}%
}}-\underbrace{S\left(  x\right)  }_{\substack{=m\left(  w\right)
-x\\\text{(by (\ref{pf.cor.id-S2.gr2.a.Sx=}))}}}\\
&  =m\left(  w\right)  -\left(  m\left(  w\right)  -x\right)  =x.
\end{align*}
In other words, $S^{2}\left(  x\right)  =x$. Hence, $\left(
\operatorname*{id}-S^{2}\right)  \left(  x\right)  =0$. Since we have shown
this for each $x\in H_{2}$, we thus conclude that $\left(  \operatorname*{id}%
-S^{2}\right)  \left(  H_{2}\right)  =0$. This proves Corollary
\ref{cor.id-S2.gr2} \textbf{(a)}.

Now we know that $\left(  \operatorname*{id}-S^{2}\right)  \left(
H_{2}\right)  =0$ (by Corollary \ref{cor.id-S2.gr2} \textbf{(a)}). In other
words, all $i\in\left\{  2,3,\ldots,2\right\}  $ satisfy $\left(
\operatorname*{id}-S^{2}\right)  \left(  H_{i}\right)  =0$ (since the only
$i\in\left\{  2,3,\ldots,2\right\}  $ is $2$). Hence, we can apply Corollary
\ref{cor.id-S2.grp} to $p=2$. Doing so, we obtain precisely the claims of
parts \textbf{(b)} and \textbf{(c)} of Corollary \ref{cor.id-S2.gr2}. (To be
precise: Corollary \ref{cor.id-S2.gr2} \textbf{(b)} follows by applying
Corollary \ref{cor.id-S2.grp} \textbf{(a)}, whereas Corollary
\ref{cor.id-S2.gr2} \textbf{(c)} follows by applying Corollary
\ref{cor.id-S2.grp} \textbf{(b)}.)
\end{proof}

\subsubsection*{Acknowledgments}

I thank Marcelo Aguiar and Amy Pang for conversations I have learnt much from.

\EditInfo{October 03, 2021}{April 21, 2022}{Pasha Zusmanovich}

\end{paper}